\documentclass[11pt]{article}

\usepackage{lineno, blindtext}
\usepackage{lipsum}
\usepackage{amsmath,amssymb,epsfig,float,color}
\usepackage{amsfonts}
\usepackage{graphicx}
\usepackage{epstopdf}
\usepackage{algorithmic}
\usepackage{xcolor}
\ifpdf
  \DeclareGraphicsExtensions{.eps,.pdf,.png,.jpg}
\else
  \DeclareGraphicsExtensions{.eps}
\fi

\title{Set title\thanks{Submitted to the editors DATE.
\funding{Funding.}}}

\author{Julio Careaga\thanks{cima 
  (\email{juliocareaga@udec.cl}).}
\and V\'ictor Osores\thanks{UCT 
  (\email{vosores@ucm.cl}).}}

\usepackage{amsopn}

\usepackage{hyperref}
\hypersetup{
    colorlinks  = true,
    linkcolor   = blue,
    filecolor   = magenta,      
    urlcolor    = cyan,
    citecolor   = red!50!black,
    }

\renewcommand{\div}{\mathrm{div}}
\newcommand{\bdiv}{\mathbf{div}}

\newcommand{\bB}[1]{\boldsymbol{#1}}
\newcommand{\mathbcal}[1]{\boldsymbol{\mathcal{#1}}}

\usepackage[]{mdframed}

\newcommand*{\bcl}{\{\mskip-6.5mu\{}
\newcommand*{\bcr}{\}\mskip-6.5mu\}}
\newcommand*{\bigbcl}{\big\{\mskip-8mu\big\{}
\newcommand*{\bigbcr}{\big\}\mskip-8mu\big\}}

\newcommand{\sgn}{\operatorname*{sgn}}

\newcommand{\upw}{{\rm Upw}}

\usepackage{bm}
\newcommand{\intvl}{\mathcal{I}}

\begin{document}

\title{A multilayer shallow water model for polydisperse reactive sedimentation\phantom{\thanks{Corresponding author}}}

\author{
{\sc Julio Careaga}\thanks{Departamento de Matem\'atica, Universidad del B\'io B\'io, Chile, email: {\tt jcareaga@ubiobio.cl}.} $^{,*}$
\quad\text{and}\quad
{\sc V\'ictor Osores}\thanks{Departamento de Matem\'atica, Universidad Cat\'olica del Maule, Chile, email: {\tt vosores@ucm.cl}.}
}

\date{ }

\maketitle
\begin{abstract}

\noindent A three-dimensional model of polydisperse reactive sedimentation is developed by means of a multilayer shallow water approach. The model consists of a variety of solid particles of different sizes and densities, and substrates diluted in water, which produce biochemical reactions while the sedimentation process occurs.
Based on the Masliyah--Lockett--Bassoon settling velocity, compressibility of the sediment and viscosity of the mixture, the system of governing equations is composed by non-homogeneous transport equations, coupled to a momentum equation describing the mass-average velocity.
Besides, the free-surface depicted by the total height of the fluid column is incorporated and fully determined through the multilayer approach. A finite volume numerical scheme on Cartesian grids is proposed to approximate the model equations. Numerical simulations of the denitrification process exemplify the performance of the numerical scheme and model under different scenarios and bottom topographies.

\end{abstract}\
\\[1ex]
\section{Introduction}

Reactive sedimentation is defined as the process in which multiple types of solid particles suspended on a fluid settle due to gravity, while reacting with substrates diluted in water. The solid phase is then composed by organic or inert matter, including bacteria considered as solid particles, and the fluid phase is constituted by substrates. An example of reactive sedimentation can be found in wastewater treatment industries, when the so-called activated sludge is treated in the Secondary Settling Tank. Activated sludge models (ASMs) have been widely developed in terms of biokinetic reactions and ordinary differential equations, see e.g.~\cite{Dupont1992, Flores2008, Gernaey2004, Henze1987WR, Henze2000ASMbook, Hu2003}. Models accounting space variations based on partial differential equations (PDEs) have been studied to a lesser extent. The one-dimensional PDE-based model presented in \cite{SDIMA_MOL}, which has evolved from the previous work \cite{SDcace_reactive}, constitutes one of the few works including the complexity of nonlinear effects appearing in the sedimentation-consolidation process. Furthermore, this model is tailored to handle mixtures of arbitrary number of solid species and liquid substrates. Extensions to Sequencing Batch Reactors, in which the main complication is related to a moving-boundary problem, can be found in \cite{SDAMM_SBR1,SDAMM_SBR2,Burger2023b}, and
comparisons with experimental data are presented in \cite{Burger2023a}.

In terms of advancements in PDE-based models of mono-dispersed and non-reactive sedimentation, including the development of suitable numerical schemes, a variety of works in different space dimensions can be found in the literature. One-dimensional approaches for batch and continuous sedimentation  such as \cite{SDcace1, Burger&E&K&L2000, Burger&K&T2005a, SDsiam3, SDwatres2, Ekama1997, Watts1996a} are mainly approximated by finite differences or finite volume schemes, with focus on non-linear flux approximations.
The special case of quasi-one-dimensional approaches \cite{Anestis1981, SDsiap_varyingA, SDcec_varyingA, Burger&D&K2004, Chancelier1994}, treat the container geometry as a coefficient function.
Already for two space dimensions, most of the models need to be supplemented with momentum balances, increasing the number of equations to solve and adding a substantial degree of difficulty.
Different variants appear depending on the velocity field considered: volume-average velocity \cite{Burger2012b,Burger&RB&T2012}, mass-average velocity \cite{Rao2002}, or solid-phase velocity \cite{Careaga2023,Gustavsson2000}. In \cite{Burger2012b}, a multiresolution finite volume scheme is used to approximate the two-dimensional model, while a combination of finite volumes and mixed finite element methods is employed in \cite{Careaga2023}. A sophisticated finite volume element method was tailored in  \cite{Burger&RB&T2012} to approximate their developed axisymmetric equations. For suspensions composed by solid species with different densities, an adequate description of the process comes from a polydisperse sedimentation description. In this line, the work in \cite{Osores2020} addresses a multilayer Saint-Venant (MSV) polydisperse model in three dimensional domains having a free-boundary.
This dynamic multilayer approach combines the theory of mixtures with the ideas presented in \cite{Audusse2005, Audusse2011, Audusse2010}, so that the final system of equations couples a type of multilayer shallow water equations with the polydisperse sedimentation model developed in \cite{Basson2009, Tory2003}.
The multilayer method provides then an efficient alternative for numerically solving the transport and flow equations by layers, reducing the problem by one space dimension. An additional gain in the formulation proposed in \cite{Osores2020} is that the vertical components of the velocity fields are computed by means of a post-processing procedure.
The numerical scheme proposed is based on the specialized methods for hyperbolic equations with non-conservative products \cite{CastroDaz2012, CastroDaz2009}, which employ the theory introduced by Dal Maso et al. \cite{DalMaso1995} to properly approximate the non-conservative products.

The aim of this work is to extend the reactive sedimentation model in \cite{SDIMA_MOL} to the three-dimensional polydisperse case following the line of \cite{Osores2020}, this is, employing an MSV approach. In contrast to the one-dimensional scenario, in which only mass balances are required, we need to take into account momentum equations. The latter contemplates the inclusion of the viscous stress tensor of the mixture in terms of a mass-average velocity.
Hence, we derive a polydisperse multidimensional reactive sedimentation model, that couples convection-diffusion-reaction transport equations with flow equations.
In turn, the MSV approach, which requires the assumption of shallow water, is employed in order to reduce the complexity that three-dimensional evolutionary PDEs mean. With this, we aspire to give a full description of the governing equations, making distinction between the vertical coordinate (pointing in the opposite direction to gravity) and the horizontal ones, so that the final equations are written in layers that subdivide the vertical axis. Moreover, we seek to develop a numerical scheme for the approximation of the model equations combining the Harten–Lax–van Leer (HLL)-type path-conservative numerical method from \cite{Castro2012} with upwind flux approximations for the substrates.

\subsection{Outline of the paper}

The paper is organized as follows. In Section~\ref{sec:modeleq}, we introduce the continuity equations for the solids and fluid phases, and momentum balance for the mixture. The governing PDEs for three-dimensional domains are written as a coupled system in terms of the solids and substrates concentrations, concentration of the mixture, mass-average velocity and pressure. In Section~\ref{sec:multilayermodel}, we develop the multilayer or shallow water version of the governing equations. We start in Section~\ref{subsec:preliminaries} introducing notation and specifications related to the layers of our multilayer approach. Then, in Section~\ref{subsec:multilayereqns}, we formulate the multilayer form of the governing equations, including a description of the normal mass fluxes. In Section~\ref{sec:numschem}, we propose a finite volume numerical scheme to approximate the multilayer version of our model. Numerical simulations of the denitrification of activated sludge process are presented in Section~\ref{sec:numexamples}. Finally, concluding remarks are summarized in Section~\ref{sec:conclusions}.

\section{Model derivation} \label{sec:modeleq}

We begin by considering a bounded domain $\Omega\subset\mathbb{R}^3$ and $\bB{x}:=(x_1,x_2,z)$ the vector of coordinates of a point in $\Omega$, where $\bB{\widetilde{x}} = (x_1,x_2)$ corresponds to its horizontal coordinates, and $\Omega_T:=\Omega\times(0,T]$ for $T>0$. The vertical direction, denoted by the $z$-coordinate, is assumed to be pointing in the opposite direction to gravity, and its unit vector is denoted by $\bB{k} = (0,0,1)$ or simply $\bB{k}=(\bB{0},1)$ with $\bB{0} = (0,0)$.
For ease in notation when describing the multiple particulate components, substrates and layers, we define for each $k\in\mathbb{N}$ the following sets on integers
$$\intvl(k):=\{1,2,\dots,k\}\quad\mbox{and} \quad\intvl_0(k):=\{0,1,\dots,k\}.$$
We consider $n_{\bB{c}}\in \mathbb{N}$ species of spherical solid particles dispersed in a viscous fluid whose vector of volume fractions is $\bB{\phi}:=\big(\phi^{(1)},\phi^{(2)},\dots,\phi^{(n_{\bB{c}})}\big)^{\tt t}\in\mathbb{R}^{n_{\bB{c}}}$, where the $i$-th phase volume fraction is denoted by $\phi^{(i)}:=\phi^{(i)}(\bB{x},t)$ for all $(\bB{x},t)\in \Omega_T$. Furthermore, we assume that each solid phase $i\in\intvl(n_{\bB{c}})$ has density $\varrho_i>0$ and particle diameter $d_i>0$, such that $d_1\geq d_2 \geq \cdots \geq d_N$. The volume fraction and density of the fluid phase is defined by $\phi_{\rm f}:=\phi_{\rm f}(\bB{x},t)$ and $\varrho_{\rm f}$ (constant) for all $(\bB{x},t)\in\Omega_T$, respectively. The velocity field of each solid phase $i\in\intvl(n_{\bB{c}})$ is denoted by
\begin{align*}
\bB{v}_i := \bB{v}_i(\bB{x},t) := (u_i(\bB{x},t),v_i(\bB{x},t),w_i(\bB{x},t))\in \mathbb{R}^3,\quad \forall (\bB{x},t)\in \Omega_T,
\end{align*}
where $\bB{\widetilde{v}}_i:=(u_i,v_i)\in \mathbb{R}^2$ corresponds to its horizontal velocity and $w_i$ is referred as the vertical component. Without loss of generality, we will also use the compact notation $\bB{v}_i = (\bB{\widetilde{v}}_i,w_i)$ to differentiate between the horizontal and vertical components of the velocity fields.
We denote the fluid velocity by $\bB{v}_{\rm f}:=\bB{v}_{\rm f}(\bB{x},t)$ and define the vector of solid concentrations by $\bB{c}:=\big(c^{(1)},c^{(2)},\dots,c^{(n_{\bB{c}})}\big)^{\tt t}\in\mathbb{R}^{n_{\bB{c}}}$ such that $c^{(i)}(\bB{x},t) := \varrho_i\phi^{(i)}(\bB{x},t)$ for all $i\in\intvl(n_{\bB{c}})$ and $(\bB{x},t)\in \Omega_T$.
The total solids volume fraction and total concentration (including the fluid phase) at each point $(\bB{x},t)\in\Omega_T$ are given by
\begin{align*}%\label{def:rho}
\phi := \phi(\bB{\phi}):=\sum_{i =1}^{n_{\bB{c}}}\phi^{(i)}=\sum_{i =1}^{n_{\bB{c}}}c^{(i)}/\rho_i,\qquad \rho:= \rho(\bB{\Phi}(\bB{x},t))  = \varrho_{\rm f}\phi_{\rm f}(\bB{x},t) + \sum_{i = 1}^{n_{\bB{c}}} \varrho_i\phi^{(i)}(\bB{x},t)\,,
\end{align*}
respectively, where $\bB{\Phi} := \big(\phi_{\rm f},\phi^{(1)},\dots,\phi^{(n_{\bB{c}})}\big)^{\tt t}$ is the vector of volume fractions including the one of the fluid phase. Using the constitutive assumption $\phi_{\rm f}+\phi = 1$, we can express $\phi_{\rm f}$ in terms of $\bB{c}$ as follows
\begin{align*}
 \phi_{\rm f}(\bB{c})= 1 - \sum_{i=1}^{n_{\bB{c}}} \dfrac{1}{\varrho_i} c^{(i)}\,,
\end{align*}
so that we can write $\bB{\Phi}$ in terms of $\bB{c}$, therefore $\bB{\Phi} = \bB{\Phi}(\bB{c})$, the same, therefore, can be done to write $\phi$ and $\rho$ in terms of $\bB{c}$. In turn, we assume that the fluid phase is composed by $n_{\bB{s}}\in\mathbb{N}$ soluble substrates diluted in water, and define the vector of substrate concentrations by $\bB{s}:=\big(s^{(1)},s^{(2)},\dots,s^{(n_{\bB{s}})}\big)^{\tt t}\in\mathbb{R}^{n_{\bB{s}}}$, where each component is a function of $(\bB{x},t)\in\Omega_T$.

Taking into account that the solid phase may react with the fluid phase (substrates), each solid specie satisfies the continuity equation or mass conservation
\begin{align}\label{eq:mass:c}
 \dfrac{\partial {c^{(i)}}}{\partial t} + \div\bigl(c^{(i)}\bB{v}_{i}\bigr) = \mathcal{R}^{(i)}_{\bB{c}}(\bB{c},\bB{s})\qquad \forall (\bB{x},t)\in\Omega_T\,,\quad\forall i\in \intvl(n_{\bB{c}})\,,
\end{align}
where $\mathcal{R}^{(i)}_{\bB{c}}:=\mathcal{R}^{(i)}_{\bB{c}}(\bB{c},\bB{s})$ is a (nonlinear) function involving the kinetic reactions between solid and fluid components. We also define the total solid reactions as $\widetilde{\mathcal{R}}_{\bB{c}} := {\mathcal{R}}_{\bB{c}}^{(1)}+{\mathcal{R}}_{\bB{c}}^{(2)}+\dots+{\mathcal{R}}_{\bB{c}}^{(n_{\bB{c}})}$. The fluid phase, on the other hand, is composed by $n_{\bB{s}}\in\mathbb{N}$ substrates and water, whose concentration is denoted by $s_{\rm w}:=s_{\rm w}(\bB{x},t)$, so that the following equation holds
\begin{align*}
 \varrho_{\rm f}\phi_{\rm f} = s^{(1)}+s^{(2)}+\dots+s^{(n_{\bB{s}})} + s_{\rm w}\,.
\end{align*}
Similarly to \eqref{eq:mass:c}, the fluid phase fulfills the mass conservation equation
\begin{align} \label{eq:mass:phif}
 \varrho_{\rm f} \dfrac{\partial \phi_{\rm f}}{\partial t} + {\rm div}(\varrho_{\rm f}\phi_{\rm f}\bB{v}_{\rm f}) = \widetilde{\mathcal{R}}_{\boldsymbol{s}}(\bB{c},\bB{s})\qquad \forall(\bB{x},t)\in\Omega_T\,,\quad\forall i\in \intvl(n_{\bB{s}})\,,
\end{align}
where the function $\widetilde{\mathcal{R}}_{\bB{s}}:=\widetilde{\mathcal{R}}_{\bB{s}}(\bB{c},\bB{s})$ corresponds to the sum of functions describing the reactions of each substrate component, defined by $\mathcal{R}_{\bB{s}}^{(i)}:=\mathcal{R}_{\bB{s}}(\bB{c},\bB{s})$ for $i \in\intvl(n_{\bB{s}})$, therefore $$\widetilde{\mathcal{R}}_{\bB{s}}: = \mathcal{R}_{\bB{s}}^{(1)}+\mathcal{R}_{\bB{s}}^{(2)}+\dots +\mathcal{R}_{\bB{s}}^{({n_{\bB{s}}})}\,.$$
In addition, assuming that the substrate components move with the same velocity $\bB{v}_{\rm f}$ as the entire fluid phase, and that the water does not react, we can establish from \eqref{eq:mass:phif} the continuity equation for each substrate component and water concentration
\begin{subequations}
\begin{alignat}{2}
&\dfrac{\partial {s^{(i)}}}{\partial t} + \div(s^{(i)}\bB{v}_{\rm f})  = \mathcal{R}_{\bB{s}}^{(i)}(\bB{c},\bB{s})\qquad && \forall (\bB{x},t)\in \Omega_T\,,\quad\forall i\in \intvl(n_{\bB{s}})\,,\label{eq:mass:s}\\[1ex]
&\dfrac{\partial s_{\rm w}}{\partial t} + {\rm div}\left(s_{\rm w}\bB{v}_{\rm f}\right)  = 0 \qquad&& \forall (\bB{x},t)\in \Omega_T\,, \label{eq:mass:w}\\[1ex]
& s^{(1)}+\dots + s^{(n_{\bB{s}})} + s_{\rm w}  = \varrho_{\rm f}\phi_{\rm f}\qquad && \forall (\bB{x},t)\in \Omega_T\,. \notag%\label{eq:sum:s}
\end{alignat}
\end{subequations}
We observe that Equation~\eqref{eq:mass:w}, which determines $s_{\rm w}$, does not need to be solved since this variable can be recovered directly using the third equation as $ s_{ \rm w} = \varrho_{\rm f}\phi_{\rm f}(\bB{c}) - \left(s^{(1)}+\dots + s^{(n_{\bB{s}})}\right)$.

In what follows, we are going to introduce constitutive assumptions over the velocity vectors in equations \eqref{eq:mass:c} and \eqref{eq:mass:s} such that we can reformulate them in terms of the slip velocities $\bB{u}_i:=\boldsymbol{v}_i-\boldsymbol{v}_{\rm f}$ for $i\in \intvl(n_{\bB{c}})$ and the so-called mass-average velocity, defined by
\begin{align*}%\label{def:mass-average}
 \bB{v} := \dfrac{1}{\rho(\bB{c})}\biggl(\varrho_{\rm f}\phi_{\rm f}\bB{v}_{\rm f} + \sum_{j=1}^{n_{\bB{s}}}c^{(j)}\bB{v}_j\biggr)\,.
\end{align*}
Straightforwardly from the definition of the slip and mass-average velocities, we can obtain the following expressions for the velocity of the solid and fluid phase
\begin{align}
\boldsymbol{v}_i & = \boldsymbol{v} + \bB{u}_i - \dfrac{1}{\rho(\bB{c})}\sum_{j=1}^{n_{\bB{c}}}  c^{(j)}\bB{u}_j \qquad \forall  i\in \intvl(n_{\bB{c}})\,, \label{eq:relmassv}\\
\bB{v}_{\rm f} & = \bB{v} - \dfrac{1}{\varrho_{\rm f}\phi_{\rm f}}\sum_{i=1}^{n_{\bB{c}}}c^{(i)}\biggl(\bB{u}_i - \dfrac{1}{\rho(\bB{c})}\sum_{j=1}^{n_{\bB{c}}}  c^{(j)}\bB{u}_j \biggr)\,. \label{eq:relslipv}
\end{align}
Next, introducing the constitutive assumption given in \cite[Equation 2.9]{Osores2020}, we can write each slip velocity as a function of the vector of solid concentrations and its gradient, i.e., $\bB{u}_{i} = \bB{u}_{i}(\bB{c},\nabla\bB{c})$ for $i\in\intvl(n_{\bB{c}})$, so that we have
\begin{align}\label{eq:sliprel}
 c^{(i)} \biggl(\bB{u}_i - \dfrac{1}{\rho(\bB{c})}\sum_{j=1}^{n_{\bB{c}}}  c^{(j)}\bB{u}_j \biggr) =
 f_i(\bB{c})\bB{k} - \bB{A}_i(\bB{c},\nabla \bB{c})\qquad \forall i\in \intvl(n_{\bB{c}})\,,
\end{align}
where $f_i(\bB{c}) := c^{(i)} v^{\rm MLB}_i(\bB{c})$ is the batch flux function of the solid component $i\in\intvl(n_{\bB{c}})$, with $v^{\rm MLB}_i(\bB{c})$ the $i$-th Masliyah--Lockett--Bassoon (MLB)  settling velocity function \cite{Masliyah1979,Lockett1979}, and $\bB{A}_{i}:=\bB{A}_{i}(\bB{c},\nabla\bB{c})$ is the compression term corresponding to the $i$-th solid specie. For each solid component $i\in \intvl(n_{\bB{c}})$, the MLB settling velocity and compression functions are given by
\begin{align*}
v^{\rm MLB}_i(\bB{c})&:= -\dfrac{gd_1^2}{18\mu_{\rm f}} v_{\rm hs}(\phi(\bB{c}))\vartheta_i(\bB{c})\,,\\[1ex] %\label{def:vLMB}\\[1ex]
   \begin{split}
 \bB{A}_i(\bB{c},\nabla \bB{c})  &:= \dfrac{d_1^2}{18\mu_{\rm f}} v_{\rm hs}(\phi(\bB{c}))  \Bigg\{  \frac{(1-\phi(\bB{c}))c^{(i)}}{\phi(\bB{c})}\big(\delta_i -
 \boldsymbol{\delta}^\mathrm{T}\bB{c}/\rho(\bB{c})\big)\nabla \sigma_\mathrm{e}(\phi(\bB{c})) \\
  & \qquad \qquad \qquad + \,\sigma_\mathrm{e}(\phi(\bB{c}))  \Bigg[  \delta_i\nabla
   \left(  \frac{c^{(i)}}{\phi(\bB{c})}\right)-c^{(i)}\sum_{j=1}^N\frac{\delta_j}{\rho(\bB{c})}\nabla\left(\frac{c^{(j)}}{\phi(\bB{c})}  \right)  \Bigg]  \Bigg\}\,,%\label{def:AA_i}\,,
   \end{split}
\end{align*}
where $\mu_{\rm f}$ is the viscosity of the pure fluid, $v_{\rm hs}:=v_{\rm hs}(\phi)$ is the hindered settling velocity, $\sigma_{\rm e}:=\sigma_{\rm e}(\phi)\geq 0$ is the effective solid stress function with compact support contained in $[0,c_{\rm max}]$. For each $i\in\intvl(n_{\bB{c}})$, the function $\vartheta_i$ is given by
\begin{align*}
\vartheta_i(\bB{c}) := \delta_i \big(\varrho_i+\phi_{\rm f}(\bB{c})\varrho_{\rm f}-c_{\rm tot}\big)  -
\sum_{k=1}^{n_{\bB{c}}}\dfrac{c^{(k)}}{\rho(\bB{c})}\delta_k \big(\varrho_k+\phi_{\rm f}(\bB{c})\varrho_{\rm f}-c_{\rm tot} \big)\,,
\end{align*}
with $c_{\rm tot} := c^{(1)}+c^{(2)}+\dots+c^{(n_{\bB{c}})}$.
Then, replacing the relation for the slip velocities \eqref{eq:sliprel} into the expressions for the solid and fluid velocities, \eqref{eq:relmassv} and \eqref{eq:relslipv}, respectively, we can finally write the fluxes corresponding to \eqref{eq:mass:c} and \eqref{eq:mass:s} in terms of $\bB{v}$, $\bB{c}$ and $\nabla\bB{c}$, as follows
\begin{alignat*}{2}
 c^{(i)}\bB{v}_i &=\, c^{(i)}\bB{v} + f_i(\bB{c})\bB{k} -\boldsymbol{A}_i(\bB{c},\nabla \bB{c})  \qquad &\forall i\in\intvl(n_{\bB{c}})\,,\\[1.5ex]
 s^{(l)}\bB{v}_{\rm f} &=\, s^{(l)}\bB{v} - \dfrac{s^{(l)}}{\varrho_{\rm f}\phi_{\rm f}(\bB{c})}\sum_{j=1}^{n_{\bB{c}}}\Bigl(f_j(\bB{c})\bB{k} - \bB{A}_j(\bB{c},\nabla \bB{c}) \Bigr)\qquad& \forall l\in\intvl(n_{\bB{s}})\,.
\end{alignat*}
Following the dimensional analysis performed in \cite{Osores2020}, for all $i\in \intvl(n_{\bB{c}})$, the horizontal components of $\bB{A}_{i}$ are neglected, so that we can replace this vector by $a^{(i)}(\bB{c},\partial_z\bB{c})\bB{k}$, where vector $\bB{a} = (a^{(i)})_{i\in\intvl(n_{\bB{c}})}$, which collects the new compression functions of each solid specie is defined as
\begin{align*}
 \bB{a}:= \bB{a}(\bB{c},\partial_z\bB{c}):= \mathbcal{D}(\bB{c}) \dfrac{\partial\bB{c}}{\partial z} \in \mathbb{R}^{n_{\bB{c}}}\,, %\label{def:A_i}
\end{align*}
and the matrix $\mathbcal{D}(\bB{c})\in \mathbb{R}^{n_{\bB{c}}\times n_{\bB{c}}}$ is determined by the coefficients
\begin{align}
\begin{aligned} \label{def:coeffD}
 [\mathcal{D}(\bB{c})]_{i,l} &:= \dfrac{d_1^2 v_{\rm hs}(\phi(\bB{c})) }{18\mu_{\rm f}\varrho_l\phi(\bB{c})} \Biggl(\big(1-\phi(\bB{c})\big)\dfrac{c^{(l)}}{\varrho_l\rho(\bB{c})}\big(\delta_l\rho(\bB{c}) - \bB{\delta}^t \bB{c}\big)\sigma_e'(\phi(\bB{c}))\\
 &\qquad  + \sigma_e(\phi(\bB{c}))\bigg(\delta_l\hat{\delta}_{il} - \dfrac{c^{(l)}\varrho_i\delta_i}{\varrho_l\rho(\bB{c})} - \frac{c^{(l)}}{\varrho_l\phi(\bB{c})\rho(\bB{c})}\big(\delta_l\rho(\bB{c})-\bB{\delta}^t \bB{c}\big)\bigg)\Bigg)\,,\quad \forall i,l\in\intvl(n_{\bB{c}})\,,
 \end{aligned}
\end{align}
with $\bB{\delta} = (\delta_1,\delta_2,\dots,\delta_{n_{\bB{c}}})^{\tt t}$, and $\hat{\delta}_{il}$ defined as the Kronecker delta, which equals to 1 if $i=l$ and zero otherwise. The total mass conservation of the mixture is obtained by adding up each component of equations \eqref{eq:mass:c} and \eqref{eq:mass:s}, and this is given by
\begin{align*}%\label{eq:mass:mixture}
 \dfrac{\partial \rho}{\partial t} + \div\bigl(\rho\bB{v}\bigr) =  \widetilde{\mathcal{R}}_{\rho}(\bB{c},\bB{s})\qquad \forall (\bB{x},t)\in \Omega_T\,,
\end{align*}
where $\widetilde{\mathcal{R}}_{\rho}(\bB{c},\bB{s}):=\widetilde{\mathcal{R}}_{\bB{c}}(\bB{c},\bB{s})+\widetilde{\mathcal{R}}_{\bB{s}}(\bB{c},\bB{s})$.
To determine the mass-average velocity $\bB{v}$, we introduce the momentum balance equation of the mixture
\begin{align*}
 \dfrac{\partial}{\partial t}(\rho\bB{v}) + \mathbf{div}\big(\rho\bB{v}\otimes\bB{v}-\mu(\phi) \bB{e}(\bB{v})\big) +\nabla p&= -\rho g \bB{k}\qquad \forall (\bB{x},t)\in\Omega_T\,,%\label{eq:momentum:mixture}
\end{align*}
where $p$ is the pressure, $\mu:=\mu(\phi)$ is the viscosity of the mixture, and $\bB{e}(\bB{v}):=\tfrac{1}{2}\big(\nabla\bB{v}+(\nabla\bB{v})^{\tt t}\big)$ is the symmetric part of the gradient. Then, the final system of equations is the following:
  Find $(\bB{c},\bB{s},(\bB{v},p))$ such that for each $\bB{x}\in\Omega$ and $t>0$
 \begin{subequations}  \label{syst:finalmodel}
 \begin{alignat}{2}
  &\dfrac{\partial \rho}{\partial t} + \div\bigl(\rho\bB{v}\bigr) =  \widetilde{\mathcal{R}}_{\rho}(\bB{c},\bB{s})\,,&\label{eq:fm:mass:mixture}\\[1ex]
  %---------------------------------------------------------
&\dfrac{\partial}{\partial t}(\rho\bB{v}) + \mathbf{div}\big(\rho\bB{v}\otimes\bB{v}-\mu(\phi) \bB{e}(\bB{v})\big) +\nabla p= -\rho g \bB{k}\,,\label{eq:fm:momentum:mixture}\\[1ex]
  %---------------------------------------------------------
 &\dfrac{\partial {c^{(i)}}}{\partial t} + \div\Bigl(c^{(i)}\bB{v} + \big(f_i(\bB{c}) - a^{(i)}(\bB{c},\partial_z\bB{c})\big)\bB{k} \Bigr) = \mathcal{R}^{(i)}_{\bB{c}}(\bB{c},\bB{s}) &\qquad \forall i\in\intvl(n_{\bB{c}})\,,\label{eq:fm:mass:c}\\[1ex]
 %%---------------------------------------------------------
 &\dfrac{\partial {s^{(l)}}}{\partial t} + \div\biggl(
%  \mathbf{f}_{\bB{s}}^{(i)}(s^{(i)},\bB{c},\nabla\bB{c},\bB{v})
 s^{(l)}\bB{v} - s^{(l)}\sum_{j=1}^{n_{\bB{c}}}\dfrac{\big(f_j(\bB{c}) - a^{(j)}(\bB{c},\partial_z \bB{c})\big)\bB{k}}{\varrho_{\rm f}\phi_{\rm f}(\bB{c})}
 \biggr) = \mathcal{R}^{(l)}_{\bB{s}}(\bB{c},\bB{s})&\qquad \forall l\in\intvl(n_{\bB{s}})\,.\label{eq:fm:mass:s}
%  %%
 \end{alignat}
 \end{subequations}
 This system is supplemented with initial conditions and zero flux boundary conditions.
 In Section~\ref{sec:multilayermodel}, we are going to introduce the multilayer approach of \eqref{syst:finalmodel}, where an extra unknown corresponding to the height of each layer is included. Nevertheless, an additional assumption on the pressure~$p$ will allow us maintaining the same number of equations as in system \eqref{syst:finalmodel}.

 We end this section by establishing some assumptions on the reaction terms. We assume that there exist constant stoichiometric matrices $\bB{\sigma}_{\bB{c}}\in \mathbb{R}^{\ell\times n_{\bB{c}}}$ and $\bB{\sigma}_{\bB{s}}\in \mathbb{R}^{\ell\times n_{\bB{s}}}$, where $\ell>0$ is the number of non-negative reaction rates. Considering that the $\ell$ reaction rates are given by the vector function $\bB{\kappa}:=\bB{\kappa}(\bB{c},\bB{s})\in \mathbb{R}^{\ell}$, the reaction terms are defined as
 \begin{align*}
  \mathbcal{R}_{\bB{c}}(\bB{c},\bB{s}) = \bB{\sigma}_{\bB{c}}\bB{\kappa}(\bB{c},\bB{s})\,,\qquad
  \mathbcal{R}_{\bB{s}}(\bB{c},\bB{s}) = \bB{\sigma}_{\bB{s}}\bB{\kappa}(\bB{c},\bB{s})\,.
 \end{align*}
In addition, to guarantee that each component of $\bB{c}$ do not exceed $c_{\max}:=\min\{\varrho_i:\,i\in\intvl(n_{\bB{c}})\}\phi_{\max}$, with $\phi_{\max}>0$ the maximum volume fraction, we require the reaction terms tending to zero as $c_{\rm tot}$ approaches $c_{\max}$: There exists an $\varepsilon > 0$ such that $\mathbcal{R}_{\bB{c}} (\bB{c},\bB{s} ) = \bB{0}$ for all $c_{\rm tot}\geq c_{\max} - \varepsilon$.

\section{Multilayer formulation} \label{sec:multilayermodel}
\subsection{Preliminaries}\label{subsec:preliminaries}
\begin{figure}[!t]
\centering
 \includegraphics[scale=1.2]{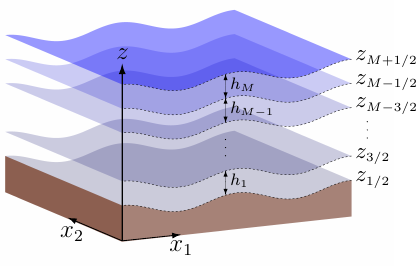}
 \caption{Illustration of the domain $\Omega$ and layers. The bottom surface (ground level) corresponds to $z_{\rm B} = z_{1/2}$ and the free-boundary is at $z_{\rm S} = z_{M+1/2}$. \label{fig:layers}}
\end{figure}

With the aim of providing a proper free-boundary description, we consider that for each $t\in[0,T]$, the domain varies with time, this is $\Omega:=\Omega(t)$.
Moreover, in order to implement a multilayer approach of \eqref{syst:finalmodel} we divide the domain $\Omega(t)$  along the vertical direction $\bB{k}$ into $M\in\mathbb{N}$ layers. More precisely, for each $\alpha \in\intvl(M)$, we consider surfaces $z:=z_{\alpha+1/2}(\widetilde{\bB{x}},t)$ with $\widetilde{\bB{x}}\in \mathbb{R}^{2}$, such that layer $\alpha$ at time $t$ corresponds to the set comprehended between $z = z_{\alpha - 1/2}(\widetilde{\bB{x}},t)$ and $z= z_{\alpha+1/2}(\widetilde{\bB{x}},t)$, i.e., $z$ lies in the interval
\begin{align*}
I_\alpha(\widetilde{\bB{x}},t):=\bigl(\,z_{\alpha-1/2}(\widetilde{\bB{x}},t),z_{\alpha+1/2}(\widetilde{\bB{x}},t)\,\bigr)\subset\mathbb{R}\,.
\end{align*}
Then, the layers and interfaces are defined by the sets
\begin{alignat*}{2}
\Omega_{\alpha}(t)&:=\Big\{(\widetilde{\bB{x}},z): \quad z\in I_\alpha(\widetilde{\bB{x}},t),\quad \widetilde{\bB{x}}\in\bB{\Pi}(t)\Big\}\,,&&\qquad \forall t\in(0,T]\,,\quad\forall\alpha\in \intvl(M)\,,\\
\Gamma_{\alpha + 1/2}(t) &: = \Big\{(\widetilde{\bB{x}},z): \quad z = z_{\alpha+1/2}(\widetilde{\bB{x}},t),\quad \widetilde{\bB{x}}\in\bB{\Pi}(t)\Big\}\,,&&\qquad \forall t\in(0,T]\,,\quad\forall\alpha\in \intvl_0(M)\,,
\end{alignat*}
where $\bB{\Pi}(t)$ maps the points in $\Omega(t)$ to their projections into the horizontal plane, i.e., $\widetilde{\bB{x}}\in \bB{\Pi}(t)$ if only if there exists $z\in \mathbb{R}$ such that $(\widetilde{\bB{x}},z)\in \Omega(t)$. The set of layers $\{\Omega_\alpha(t)\}_{\alpha\in\intvl(M)}$ fulfills the following identities
\begin{align*}
 \overline{\Omega(t)} = \bigcup_{\alpha=1}^M \overline{\Omega_{\alpha}(t)}\,,\quad \partial \Omega_{\alpha}(t) = \Gamma_{\alpha - 1/2}(t) \cup \Gamma_{\alpha + 1/2}(t) \cup \Big\{(\widetilde{\bB{x}},z): \quad \widetilde{\bB{x}}\in\partial\bB{\Pi}(t),\quad z\in I_\alpha(\widetilde{\bB{x}},t)\Big\}\,.
\end{align*}
In addition, we assume that the $M$ interfaces are of class $C^1$ in time and space.
The thickness of layer $\alpha$ is defined by the function $h_{\alpha}:=h_{\alpha}(\widetilde{\bB{x}},t) := z_{\alpha+1/2}(\widetilde{\bB{x}},t)-z_{\alpha - 1/2}(\widetilde{\bB{x}},t)$ which will vary with respect to the horizontal position $\widetilde{\bB{x}}$ and time $t$, see Figure~\ref{fig:layers}. The top and bottom surfaces are denoted by $z_{\rm B}:=z_{1/2}$ and $z_{\rm S} := z_{M+1/2}$, and the height of the fluid is $h:=z_{\rm S}-z_{\rm B}$, and there hold
\begin{align*}
 h= \sum_{\alpha=1}^{M}h_{\alpha}\,,\quad z_{M+1/2} & = z_{\rm B} + h\quad\mbox{and}\quad z_{\alpha+1/2} = z_{\rm B}+ \sum_{j=1}^{\alpha}h_j\,,\quad \forall \alpha=\intvl(M)\,.
\end{align*}
Furthermore, we assume that the layer thicknesses are small enough to neglect the dependence of the horizontal velocities and the concentrations on the vertical variable inside each layer. Finally, for each $\alpha\in \intvl(M)$ and function $\varphi:\Omega\to \mathbb{R}$, we define the one-sided limits
\begin{align*}
 \varphi_{\alpha+ 1/2}^- := (\varphi|_{\Omega_{\alpha}})|_{\Gamma_{\alpha+1/2}},\qquad 
 \varphi_{\alpha+ 1/2}^+ := (\varphi|_{\Omega_{\alpha+1}})|_{\Gamma_{\alpha+1/2}}\,.
\end{align*}
If $\varphi$ is continuous across $\Gamma_{\alpha + 1/2}(t)$, we simply set $\varphi_{\alpha + 1/2}:= \varphi|_{\Gamma_{\alpha+1/2}}$.\\

\subsection{Multilayer model equations}\label{subsec:multilayereqns}

In what follows, we are going to derive a multilayer approach of system~\eqref{syst:finalmodel}, so that the balance equations are established on each layer $\Omega_{\alpha}(t)$ for all $\alpha\in\intvl(M)$. In this regard, to properly write every equation in \eqref{syst:finalmodel} over each $\alpha$-layer, we need to take into account the flux transmissions at the interfaces $\Gamma_{\alpha+1/2}$. Before delving deeper, we introduce the divergence and gradient operators in the horizontal coordinates as $\div_{\bB{x}}$ and $\nabla_{\bB{x}}$, respectively, this is
\begin{align*}
 \div_{\bB{x}} := \dfrac{\partial}{\partial x_1}+\dfrac{\partial}{\partial x_2}\,,\qquad \nabla_{\bB{x}}:=(\partial_{x_1},\partial_{x_{2}})^{\tt t}\,,
\end{align*}
and define $\bdiv_{\bB{x}}$ as the row-wise tensorial version of the divergence $\div_{\bB{x}}$. 
For all $\alpha \in \intvl(M)$, we employ the subscript $\alpha$ to denote the restriction of each function to the layer $\alpha$, i.e., $\bB{c}_\alpha=\bB{c}|_{\Omega_{\alpha}(t)}$, $\bB{s}_\alpha=\bB{s}|_{\Omega_{\alpha}(t)}$ and $\bB{v}_{\alpha} = \bB{v}|_{\Omega_{\alpha}(t)}$, and similarly with the rest of scalar, vector and tensor variables. In addition, we split the volume-average velocity at the $\alpha$-layer  into its horizontal component $\bB{\widetilde{v}}_{\alpha}$ and vertical velocity $w_{\alpha}$, this is $\bB{v}_{\alpha} := (\bB{\widetilde{v}}_{\alpha},w_\alpha)$. Among the assumptions required by this approach, we assume the layer thicknesses are small enough to neglect the dependence of the horizontal velocities, and the concentration of each species on the vertical coordinate within each layer. This means that
\begin{align*}%\label{u_constant}
\partial_z \bB{\widetilde{v}}_{\alpha} = \bB{0}\quad \mbox{and} \quad \partial_z \bB{c}_{\alpha} =\partial_z \bB{s}_{\alpha}= \bB{0}\qquad \forall \alpha \in \intvl(M)\,.
\end{align*}
In addition, we assume that the vertical volume-average velocity $w_\alpha$ is piecewise linear in $z$, and possibly discontinuous. Under this assumption the vertical velocity and (hydrostatic) pressure are piecewise linear in~$z$, i.e., 
\begin{align*}%\label{strucvp}
\begin{split} 
\partial_z w_{\alpha} & = (\partial_z w_{\alpha}) (\bB{x},t)   \quad \mbox{and}\quad \partial_{z} p_{\alpha} = (\partial_{z} p_{\alpha})(\bB{x},t)\qquad \forall \alpha\in \intvl(M)\,.
 \end{split}   
\end{align*}
More precisely, the assumption of a hydrostatic pressure means that
\begin{align*}
 p_{\alpha}(\bB{x},z,t) = p_{\alpha+1/2} (\bB{x},t) + \rho_{\alpha}g(z_{\alpha+1/2}-z),\qquad p_{\alpha+1/2}(\bB{x},t) := p_{\rm S}(\bB{x},t) + g\sum_{\beta=\alpha+1}^{M}\rho_\beta h_{\beta}(\bB{x},t)\,,
\end{align*}
with $p_{\rm S}$ the pressure at the free surface.

We begin by describing the multilayer version of equation \eqref{eq:fm:mass:mixture} for the total density of the mixture $\rho_\alpha$ at the layer $\alpha\in\intvl(M)$. Integrating \eqref{eq:fm:mass:mixture} with respect to the vertical coordinate over the interval $I_{\alpha}$, we obtain
\begin{align}\label{eq:PDErho}
 \partial_{t}(\rho_{\alpha} h_{\alpha}) + \div_{\bB{x}}(\rho_{\alpha} h_{\alpha}\bB{\widetilde{v}}_{\alpha}) = \mathcal{G}_{\alpha+1/2}^{-}-\mathcal{G}_{\alpha-1/2}^{+}+
 h_{\alpha}
 \big(\widetilde{\mathcal{R}}_{\bB{c}}(\bB{c}_{\alpha},\bB{s}_{\alpha})+\widetilde{\mathcal{R}}_{\bB{s}}(\bB{c}_{\alpha},\bB{s}_{\alpha})\big)\,,
\end{align}
where the normal mass fluxes across the interfaces, defined by $\mathcal{G}_{\alpha+1/2}$, are assumed to satisfy the  continuity condition $\mathcal{G}_{\alpha+1/2}^{-} = \mathcal{G}_{\alpha+1/2}^{+} = \mathcal{G}_{\alpha+1/2}$ for each $\alpha \in \intvl(M)$, with
\begin{align}\label{def:flux:G}
 \mathcal{G}_{\alpha+1/2} := \rho_{\alpha}\big(\partial_tz_{\alpha+1/2} + \bB{\widetilde{v}}_{\alpha}^{\tt t}\cdot \nabla_{\bB{x}}z_{\alpha+1/2} - w_{\alpha}\big)\,.
\end{align}
An explicit expression for this flux without including $w_\alpha$ will be provided later on in this section.
Similarly, the multilayer version of equation \eqref{eq:fm:mass:c}, which gives the vector of concentrations of solid species $\bB{c}_{\alpha}$ at the layer $\alpha\in \intvl(M)$, is obtained by integration on the vertical axis
\begin{align}\label{eq:PDEc}
 \partial_{t}(c^{(i)}_{\alpha} h_{\alpha}) + \div_{\bB{x}}(c^{(i)}_{\alpha} h_{\alpha}\bB{\widetilde{v}}_{\alpha}) = \mathcal{G}_{\alpha+1/2}^{(i)}-\mathcal{G}_{\alpha-1/2}^{(i)}+h_{\alpha}\mathcal{R}^{(i)}_{\bB{c}}(\bB{c}_{\alpha},\bB{s}_\alpha)\,,\quad \forall i\in \intvl(n_{\bB{c}})\,,
\end{align}
where the normal mass fluxes of each $i$-th solid specie, which are denoted by $\smash{\mathcal{G}_{\alpha+1/2}^{(i)}}$, are defined in the same way as in \cite{Osores2019}, such that they are continuous across the interface $\Gamma_{\alpha+1/2}$ and
\begin{align*}
 \sum_{i=1}^{n_{\bB{c}}} \mathcal{G}_{\alpha+1/2}^{(i)} = \mathcal{G}_{\alpha+1/2}\,.
\end{align*}
Then, the following relation between the normal mass fluxes is derived
\begin{align}
 \mathcal{G}_{\alpha+1/2}^{(i)} =  \tilde{c}^{(i)}_{\alpha+1/2}\mathcal{G}_{\alpha+1/2} - \tilde{f}^{(i)}_{\alpha+1/2} + \tilde{a}^{(i)}_{\alpha+1/2}\qquad \forall i\in \intvl(n_{\bB{c}}),\quad \forall \alpha\in \intvl_0(M)\,, \label{def:vfluxGi}
\end{align}
where the averages $\tilde{c}^{(i)}_{\alpha+1/2}$, $\tilde{f}^{(i)}_{\alpha+1/2}$, and $\tilde{a}^{(i)}_{\alpha+1/2}$ are defined by
\begin{align*}
 & \tilde{c}^{(i)}_{\alpha+1/2}:= \frac{1}{2}\left(\frac{c^{(i)}_{\alpha+1}}{\rho_{\alpha +1}}+\frac{c^{(i)}_{\alpha}}{\rho_\alpha} \right)\,,\quad
 \tilde{f}^{(i)}_{\alpha+1/2}:= \frac{1}{2}\Big( f^{(i),+}_{\alpha+1/2} + f^{(i),-}_{\alpha+1/2}\Big)\,, \quad
\tilde{a}^{(i)}_{\alpha+1/2}:= \frac{1}{2}\Big( a^{(i),+}_{\alpha+1/2} + a^{(i),-}_{\alpha+1/2}\Big)\,.
\end{align*}
Employing vector notation, we define the vector of normal mass fluxes for the solid components as $\bB{\mathcal{G}}_{\alpha+1/2}:=\big(\mathcal{G}_{\alpha+1/2}^{(i)}\big)_{i\in\intvl(n_{\bB{c}})}$ and the vectors of the average variables as
\begin{align*}
\bB{\tilde{c}}_{\alpha+1/2}:=\big(\tilde{c}^{(i)}_{\alpha+1/2}\big)_{i\in\intvl(n_{\bB{c}})}\,,\qquad \bB{\tilde{f}}_{\alpha+1/2}:=\big(\tilde{f}^{(i)}_{\alpha+1/2}\big)_{i\in\intvl(n_{\bB{c}})}\,,\qquad
\bB{\tilde{a}}_{\alpha+1/2}:=\big(\tilde{a}^{(i)}_{\alpha+1/2}\big)_{i\in\intvl(n_{\bB{c}})}\,.
\end{align*}
The multilayer approach of equation \eqref{eq:fm:mass:s} given for the vector of substrate concentrations $\bB{s}_{\alpha}$ at each layer $\alpha\in\intvl(M)$, is obtained performing similar calculations as before, this is
\begin{align}\label{eq:PDEs}
 \partial_{t}(s^{(l)}_{\alpha} h_{\alpha}) + \div_{\bB{x}}(s^{(l)}_{\alpha} h_{\alpha}\bB{\widetilde{v}}_{\alpha}) = \mathcal{S}_{\alpha+1/2}^{(l)}-\mathcal{S}_{\alpha-1/2}^{(l)}+h_{\alpha}\mathcal{R}^{(l)}_{\bB{s}}(\bB{c}_{\alpha},\bB{s}_\alpha) \qquad \forall l\in\intvl(n_{\bB{s}})\,,
\end{align}
where the normal mass flux corresponding to the $l$-substrate is denoted by $\mathcal{S}_{\alpha+1/2}^{(l)}$, which is defined as
\begin{align}
 \mathcal{S}_{\alpha+1/2}^{(l)} :=
 \tilde{s}_{\alpha+1/2}^{(l)}\mathcal{G}_{\alpha+1/2} +
 \dfrac{1}{2\varrho_{\rm f}}\bigg(
 \dfrac{s^{(l)}_{\alpha}}{\phi_{{\rm f},\alpha}}+
 \dfrac{s^{(l)}_{\alpha+1}}{\phi_{{\rm f},\alpha+1}}
 \bigg)
 \sum_{j=1}^{n_{\bB{c}}} \big(\tilde{f}^{(j)}_{\alpha+1/2} -
 \tilde{a}^{(j)}_{3,\alpha+1/2}\big)\,, \label{def:vfluxSi}%\qquad \forall i\in \intvl(n_{\bB{s}})\,.
\end{align}
with $\tilde{s}_{\alpha+1/2}^{(l)}:= \tfrac{1}{2}\big(s^{(l)}_{\alpha}/\rho_{\alpha}+s^{(l)}_{\alpha+1}/\rho_{\alpha+1}\big)$. In the same way as we did for the solid phase, we introduce the vectors $\bB{\tilde{s}}_{\alpha+1/2} := \big(\tilde{s}_{\alpha+1/2}^{(l)}\big)_{l\in\intvl(n_{\bB{s}})}$ and $\bB{\mathcal{S}}_{\alpha+1/2}:=\big(\mathcal{S}^{(i)}_{\alpha+1/2}\big)_{i\in\intvl(n_{\bB{s}})}$.
For the multilayer version of the momentum equation \eqref{eq:fm:momentum:mixture}, we observe that such a formulation can be obtained in the same way as in \cite{Osores2019}, so that we skip its derivation and directly introduce it as follows 
\begin{align}\label{eq:PDEv}
 \begin{aligned}
 &\partial_{t} (h_{\alpha} \rho_{\alpha} \bB{\widetilde{v}}_{\alpha} )
   +   \nabla_{\bB{x}}   \cdot   \bigl(h_{\alpha} \rho_{\alpha}
  \bB{\widetilde{v}}_{\alpha} \otimes \bB{\widetilde{v}}_{\alpha} \bigr)
     +  h_\alpha\bigl( \nabla_{\bB{x}}  \bar{p}_\alpha+ g \rho_\alpha 
    \nabla_{\bB{x}} \bar{z}_\alpha\bigr) \\ 
    &\qquad =  \dfrac{1}{2} \mathcal{G}_{\alpha+1/2}\big(
     \bB{\widetilde{v}}_{\alpha+1}-\bB{\widetilde{v}}_\alpha\big) - \dfrac{1}{2}\mathcal{G}_{\alpha-1/2}\big(\bB{\widetilde{v}}_{\alpha}-\bB{\widetilde{v}}_{\alpha-1}\big)
 + \bB{K}_{\alpha+1/2}-\bB{K}_{\alpha-1/2}\,,
 \end{aligned}
\end{align}
where the vector fluxes $\bB{K}_{\alpha+1/2}$ for $\alpha\in\intvl_0(M)$, which appears after neglecting terms on the viscous stress tensor, is given by
\begin{align}\label{def:Kflux}
 \bB{K}_{\alpha+1/2}:= \dfrac{\mu(\phi_{\alpha+1/2})}{2}  \frac{\widetilde{\bB{v}}_{\alpha+1}-\widetilde{\bB{v}}_\alpha}{h_{\alpha+1/2}}\,,
\end{align}
with $h_{\alpha+1/2}$ is the distance between the midpoints of layers $\alpha$ and $\alpha+1/2$. On the other hand $\bar{p}_{\alpha}$ and $\bar{z}_{\alpha}$ are defined by \eqref{def:pbar} and \eqref{def:zbar}, respectively.

Next, in order to close the system, we further assume that the thickness of each layer~$h_{\alpha}$ is a fixed fraction of the total height~$h$, such that  $h_{\alpha} = l_{\alpha} h$  for $\alpha \in \intvl(M)$, where $l_1, \dots, l_M$ are positive numbers that add up to 1. Then, defining new variables $m_\alpha:=\rho_\alpha h$, $\bB{q}_{\alpha} := \rho_{\alpha} h
\bB{\widetilde{v}}_{\alpha}$, $\bB{r}_{\alpha} :=\bB{c}_{\alpha} h$ and $\bB{\zeta}_{\alpha} :=\bB{s}_{\alpha} h$ for $\alpha \in \intvl(M)$ and replacing them into \eqref{eq:PDErho}, \eqref{eq:PDEv}, \eqref{eq:PDEc} and \eqref{eq:PDEs}, we can write the multilayer formulation of \eqref{syst:finalmodel} as: For all $\alpha\in\intvl$, find $(m_{\alpha},\bB{q}_{\alpha},\bB{r}_{\alpha},\bB{\zeta}_{\alpha})$ such that for each $\bB{\widetilde{x}}\in\bB{\Pi}(t)$ and $t>0$
\begin{subequations} \label{syst:finalmodel_ml} 
\begin{align}
&\partial_{t} m_\alpha + \div_{\bB{x}}(\bB{q}_{\alpha}) = \frac{1}{l_\alpha}\big(\mathcal{G}_{\alpha+1/2} - \mathcal{G}_{\alpha-1/2}\big) + 
h\widetilde{\mathcal{R}}_{\rho,\alpha}\,,\label{eq:ml:m}\\[1ex]
& \begin{array}{l}
\displaystyle\partial_{t} \bB{q}_{\alpha}
   +   \div_{\bB{x}} \bigg(\frac{\bB{q}_{\alpha} \otimes \bB{q}_{\alpha}}
  {m_\alpha}\bigg)
     +  h\bigl( \nabla_{\bB{x}}  \bar{p}_\alpha+ g \rho_\alpha 
    \nabla_{\bB{x}} \bar{z}_\alpha\bigr)\\[1ex]
\displaystyle\qquad
= \frac{1}{l_\alpha}\big(\bB{\tilde{q}}_{\alpha+1/2} \mathcal{G}_{\alpha+1/2}- \bB{\tilde{q}}_{\alpha-1/2} \mathcal{G}_{\alpha-1/2}) + \frac{1}{l_\alpha}\big(\bB{K}_{\alpha+1/2}-\bB{K}_{\alpha-1/2}\big)\,,\label{eq:ml:v}
\end{array}\\[1ex]
&\partial_{t}r^{(i)}_{\alpha} + \div_{\bB{x}}\bigg(\frac{r^{(i)}_{\alpha}\bB{q}_\alpha}{m_\alpha}\bigg) =
\frac{1}{l_\alpha}\big(\mathcal{G}_{\alpha+1/2}^{(i)} - \mathcal{G}_{\alpha-1/2}^{(i)}\big) + 
h\mathcal{R}^{(i)}_{\bB{c},\alpha} \qquad \forall i \in \intvl(n_{\bB{c}})\,,\label{eq:ml:c}\\[1ex]
&\partial_{t}\zeta^{(i)}_{\alpha} + \div_{\bB{x}}\bigg(\frac{\zeta^{(i)}_{\alpha}\bB{q}_\alpha}{m_\alpha}\bigg) =
\frac{1}{l_\alpha}\big(\mathcal{S}_{\alpha+1/2}^{(i)} - \mathcal{S}_{\alpha-1/2}^{(i)}\big) + 
h\mathcal{R}^{(i)}_{\bB{s},\alpha} \qquad \forall i \in \intvl(n_{\bB{s}})\,,\label{eq:ml:s}
\end{align} 
\end{subequations}
where $\bar{p}_\alpha$, $\bar{z}_{\alpha}$, and average $\bB{\tilde{q}}_{\alpha+1/2}$ are defined by
\begin{align}
&\bar{p}_\alpha := p_{\mathrm{S}} + \frac{1}{2}gl_\alpha m_\alpha + \sum_{\beta=\alpha+1}^M gl_\beta m_\beta\,,\label{def:pbar}\\
&\bar{z}_\alpha  := z_\mathrm{B} +\dfrac{1}{2} l_\alpha h + \sum_{\beta=1}^{\alpha-1} l_\beta h\,,\label{def:zbar}\\
&\bB{\tilde{q}}_{\alpha +1/2} :=\frac{1}{2} \biggl( 
 \frac{\bB{q}_{\alpha +1}}{m_{\alpha +1}}+\frac{\bB{q}_{\alpha}}{m_\alpha}\biggr)\,. \notag%\label{def:vbar}
\end{align}
Using definitions \eqref{def:pbar} and \eqref{def:zbar}, together with $l_1+\dots+l_M=1$ and $m_{\alpha} = \rho_{\alpha}h$, the third term on the left-hand side of \eqref{eq:ml:v}, accounting for the pressure gradient, can be written as
\begin{gather*}
\begin{split}%\label{def:psi}
 \bB{\psi}_{\alpha} := h\bigl( \nabla_{\bB{x}}  \bar{p}_\alpha+ g \rho_\alpha
    \nabla_{\bB{x}} \bar{z}_\alpha\bigr) & = g m_\alpha\nabla_{\bB{x}}(z_\mathrm{B}+h) + gh^2\biggl(\frac{1}{2}l_\alpha + \sum_{\beta=\alpha+1}^M l_\beta \biggr)\nabla_{\bB{x}}\rho_{\alpha}\\
    &\qquad + gh\sum_{\beta=\alpha+1}^M l_\beta  \nabla_{\bB{x}}\big(m_\beta-m_{\alpha}\big)\,.
    \end{split}
\end{gather*}
The conservation equation for the total mass of the solid phase is obtained by adding up equation \eqref{eq:ml:m} (after multiplying by $l_{\alpha}$ in both sides of the equation) from $\alpha=0$ to $\alpha=M$, this is 
\begin{align} \label{eq:ml:mtot}
\partial_t \bar{m}  + \div_{\bB{x}}\Biggl(\sum_{\beta=1}^M l_\beta \bB{q}_{\beta}\Biggr)
 =\mathcal{G}_{M+1/2}-\mathcal{G}_{1/2} + \sum_{\beta=1}^{M}l_{\beta}h\widetilde{\mathcal{R}}_{\rho,\beta} \,,
\end{align}
where $\smash{\bar{m}:=\sum_{\beta=1}^M l_\beta m_\beta}$ is the total mass of the solid phase, and~$\mathcal{G}_{M+1/2}$ and~$\mathcal{G}_{1/2}$ represent the mass transfer on the free surface and at the bottom, respectively. We observe that the mass at each layer and height can be written in terms of $\bar{m}$ and the concentration of solid species as follows:
\begin{align*}
 m_{\alpha} & = \bar{m}
 + \sum_{\beta=1}^{M}\sum_{l=1}^{n_{\bB{c}}} \dfrac{\varrho_l-\varrho_{\rm f}}{\varrho_l} l_{\beta}\big(r^{(l)}_{\alpha}-r^{(l)}_{\beta}\big)\,,\\[1ex]% \label{eq:formula:m}
 h & = \dfrac{1}{\varrho_{\rm f}}\biggl(\bar{m}
 - \sum_{\beta=1}^{M}\sum_{l=1}^{n_{\bB{c}}} \dfrac{\varrho_l-\varrho_{\rm f}}{\varrho_l} l_{\beta}r^{(l)}_{\beta}\biggr)\,. %\label{eq:formula:h}
\end{align*}
For each $\alpha\in\intvl(M)$, the normal mass flux $\mathcal{G}_{\alpha+1/2}$, as it was defined in \eqref{def:flux:G}, depends explicitly on the vertical velocities $w_{\alpha}$, which we are going to be determined on a post processing procedure. To avoid this dependency, we reformulate these fluxes in the same way as in \cite{Osores2020}. First, we define the following vector and scalar variables
\begin{align}\nonumber
&\boldsymbol{R}_{\alpha}:=\bB{q}_{\alpha}-\sum_{j=1}^{n_{\bB{c}}} r_{\alpha}^{(j)}\frac{\bB{q}_{\alpha}}{m_{\alpha}} \frac{\varrho_j-\varrho_{\rm f}}{\varrho_j}, \quad
 \bar{\boldsymbol{R}}:=\sum_{\beta=1}^Ml_\beta \boldsymbol{R}_{\beta},\quad
 \tilde{\rho}_{\alpha+1/2} := \dfrac{2\rho_{\alpha}\rho_{\alpha+1}}{\rho_{\alpha}+\rho_{\alpha+1}},\quad 
 L_{\alpha} := \sum_{\beta=1}^{\alpha} h_{\beta}\,,
\end{align}
for all $\alpha\in \intvl(M)$.  
Then, following the same calculations as in \cite{Osores2020} we obtain
\begin{align} \label{eq:fluxGnew}
\mathcal{G}_{\alpha+1/2} = \frac{\tilde{\rho}_{\alpha+1/2}}{\varrho_{\rm f}}  \sum_{\beta=1}^{\alpha} l_{\beta} \div_{\bB{x}} (\boldsymbol{R}_{\beta}-\bar{\boldsymbol{R}})+ \tilde{\mathcal{G}}_{\alpha+1/2}^{\mathcal{R}} + \tilde{\mathcal{G}}_{\alpha+1/2}   \qquad \forall \alpha \in \intvl_0(M)\,,
\end{align}
where $\tilde{\mathcal{G}}_{\alpha+1/2}^{\mathcal{R}}$ corresponds to the contribution made by the reaction terms, given by
\begin{align*}
\tilde{\mathcal{G}}_{\alpha+1/2}^{\mathcal{R}}&:=
\sum_{j=1}^{n_{\bB{c}}}\dfrac{(\varrho_j-\varrho_{\rm f})\tilde{\rho}_{\alpha+1/2}}{\varrho_j\varrho_{\rm f}}\bigg(\sum_{\beta=1}^{\alpha} l_{\beta}h \mathcal{R}_{\bB{c},\beta}^{(j)}
- L_{\alpha}\sum_{\gamma=1}^{M} l_{\gamma}h \mathcal{R}_{\bB{c},\gamma}^{(j)}\bigg)\\
&\qquad - \dfrac{\tilde{\rho}_{\alpha+1/2}}{\varrho_{\rm f}} \bigg( \sum_{\beta=1}^{\alpha}l_{\beta}h\widetilde{\mathcal{R}}_{m,\beta}
- L_{\alpha}\sum_{\gamma=1}^{M}l_{\gamma}h\widetilde{\mathcal{R}}_{\rho,\gamma}\bigg)\,.
\end{align*}
We assume that the reaction terms affect only locally on each layer, so that we neglected this vertical contribution and simply set $\tilde{\mathcal{G}}_{\alpha+1/2}^{\mathcal{R}} \equiv 0$ for all $\alpha\in\intvl_0(M)$. The third term in \eqref{eq:fluxGnew}, which is related to the nonlinear hindered-settling fluxes and compression functions, is defined as
\begin{align}
\begin{aligned}\label{def:vfluxGtilde}
\tilde{\mathcal{G}}_{\alpha+1/2} &
:=\sum_{j=1}^{n_{\bB{c}}} \dfrac{(\varrho_j-\varrho_{\rm f})\tilde{\rho}_{\alpha+1/2}}{\varrho_j\varrho_{\rm f}}
\biggl( -\left(\tilde{f}_{\alpha+1/2}^{(j)}-\tilde{a}_{\alpha+1/2}^{(j)}\right)
+ (1-L_{\alpha})\left( \tilde{f}_{1/2}^{(j)}-\tilde{a}_{1/2}^{(j)}\right) \\
& \qquad + L_{\alpha}\left( \tilde{f}_{M+1/2}^{(j)}-\tilde{a}_{M+1/2}^{(j)}\right)
\biggr) + (1-L_{\alpha})\frac{\tilde{\rho}_{\alpha+1/2}}{\tilde{\rho}_{1/2}}\mathcal{G}_{1/2}
           + L_{\alpha}\frac{\tilde{\rho}_{\alpha+1/2}}{\tilde{\rho}_{M+1/2}}\mathcal{G}_{M+1/2}\,.
\end{aligned}
\end{align}
In the special case of zero flux boundary conditions at $z=z_B$ and $z=h$, we have $\mathcal{G}_{1/2} = \mathcal{G}_{M+1/2} = 0$, and the formula above reduces to
\begin{align*}
 \tilde{\mathcal{G}}_{\alpha+1/2} &
:= - \sum_{j=1}^{n_{\bB{c}}} \dfrac{(\varrho_j-\varrho_{\rm f})\tilde{\rho}_{\alpha+1/2}}{\varrho_j\varrho_{\rm f}}
\left(\tilde{f}_{\alpha+1/2}^{(j)}-\tilde{a}_{\alpha+1/2}^{(j)}\right)\,.
\end{align*}

\noindent \emph{Computation of the vertical velocity}. The multilayer formulation presented in this section does not involve the vertical velocity of the mixture $w$. Nevertheless, this variable can be computed by means of a recursive post-procesing procedure, as the one presented in \cite{Osores2020}. In this way, the procedure reads as follows. The first step is to recursively calculate the vertical velocities at the interfaces for all $\alpha\in\intvl(M)$, this is
\begin{align*}
 w_{1/2}^{+} & = \partial_{t}z_{\mathrm{B}}+\bB{\widetilde{v}}_{1}^{\tt t}\cdot\nabla_{\bB{x}}z_{\mathrm{B}} - \frac{G_{1/2}}{\rho_1}\,,\\
 %-----------------
 w_{\alpha+1/2}^{-} & = w_{\alpha-1/2}^{+} - \frac{h_\alpha}{\rho_{\alpha}} \bigl( \partial_{t} \rho_{\alpha} + \div_{\bB{x}}(\rho_{\alpha}\bB{\widetilde{v}}_{\alpha}) - \widetilde{\mathcal{R}}_{\rho,\alpha}\bigr)\,,\\
 %-----------------
 w_{\alpha+1/2}^{+} & =
 \frac{1}{\rho_{\alpha}}\Bigl(
 (\rho_{\alpha}-\rho_{\alpha-1})\partial_{t}z_{\alpha-1/2} +
  (\rho_{\alpha+1}\bB{\widetilde{v}}_{\alpha+1}-\rho_{\alpha}\bB{\widetilde{v}}_{\alpha})^{\tt t}
  \cdot\nabla_{\bB{x}}z_{\alpha-1/2} + \rho_{\alpha -1} w_{\alpha-1/2}^{-}\Bigr)\,.
\end{align*}
Then, for $\alpha \in\intvl(M)$ and $z \in (z_{\alpha-1/2}, z_{\alpha+1/2})$ we compute the vertical velocities by:
\begin{align*}%\label{v:2}
\begin{split}
 w_{\alpha}(\bB{x},z,t)& = w_{\alpha-1/2}^+-\frac{1}{\rho_{\alpha}}  \bigl( \partial_{t} \rho_{\alpha} + \div_{\bB{x}}(\rho_{\alpha}\bB{\widetilde{v}}_{\alpha}) - \widetilde{\mathcal{R}}_{\rho,\alpha} \bigr)
 (z-z_{\alpha-1/2})\,.
\end{split}
\end{align*}

\section{Numerical scheme}\label{sec:numschem}

In this section, we describe the finite volume numerical scheme designed to approximate system \eqref{syst:finalmodel_ml} for the case of $\Omega\subset\mathbb{R}^3$ and two-dimensional horizontal domain projections. We employ the numerical scheme described in \cite{Osores2019} to our model, which uses the HLL-PVM-1U method developed in \cite{Castro2012} combined with upwind flux approximations. Besides, the viscous terms in \eqref{eq:ml:v} are computed through a splitting procedure.

In a Cartesian grid, we consider $N$ control volumes $V_{i}\subset\mathbb{R}^2$ for $i\in\intvl(N)$, given by squares of size $\Delta x\times\Delta y$ with $\Delta x$ and $\Delta y$ positive, such that $\mathcal{T}:=\smash{\{V_i\}_{i=1}^{N}}$ is a partition of the horizontal domain. We assume that the horizontal domain remains constant, so that the partition $\mathcal{T}$ remains constant over time. In addition, we define $e_{i,j}$ as the edge between two adjacent control volume~$V_i$ and~$V_j$, and $\bB{\eta}_{i,j}$~as the unitary normal vector at $e_{i,j}$ pointing from~$V_i$ to~$V_j$. In turn, the set of indexes of control volumes neighboring $V_i$ is defined by~$J_i$. The center of mass of the $i$-th control volume is denoted by~$\bB{x}_i$ while $|V_i|$ and~$|e_{i,j}|$ stand for the area of~$V_i$ and length of~$e_{i,j}$, respectively. Given a function $\varphi$, we define $\bcl\varphi\bcr_{i,j}:= (\varphi_i+\varphi_j)/2$, the average of $\varphi$ at the edge $e_{i,j}$, which is defined analogously for vector and tensor variables. For a time span specified by $t_n = n\Delta t$ with $\Delta t>0$ being a time step and $n\in \mathbb{N}$, we approximate each unknown by its volume average. Therefore, $m_{\alpha}$ is approximated at the $i$-th control volume $V_i$ and time $t = t_n$ by
\begin{align*}
 m_{i,\alpha}^n = \dfrac{1}{|V_i|}\int_{V_i} m_{\alpha}(\cdot,t_n)\,{\rm d}\bB{x}\qquad \forall \alpha\in \intvl(M)\,,
\end{align*}
and analogously with the rest of unknowns. We start by describing the approximation of the horizontal fluxes of the system, which are related to the horizontal divergence and gradient operators in \eqref{syst:finalmodel_ml}. For the first two equations \eqref{eq:ml:m} and \eqref{eq:ml:v}, given $i\in\intvl(N)$ and $j\in J_i$, we approximate their respective horizontal fluxes on the interface $e_{i,j}$ by
\begin{align*}
 \mathcal{F}_{i,j,\alpha}^{m,n} & =
 \bigbcl\bB{q}_{\alpha}^n\bigbcr_{i,j}^{\tt t}
\cdot\bB{\eta}_{i,j}
- \frac{1}{2} \Bigl( \theta^n_{0,i,j}\big( m^n_{j,\alpha}
- m^n_{i,\alpha} + b^n_{i,j,\alpha} \big)
+ \theta^n_{1,i,j}  (\bB{q}^n_{j,\alpha} - \bB{q}^n_{i,\alpha})^{\tt t}\cdot\bB{\eta}_{i,j}  \Bigr)\,,\\[1ex]
%--------------------------------------------
\mathbcal{F}^{\bB{q},n}_{i,j,\alpha} & = \bigbcl m_\alpha^{-1}\bB{q}_{\alpha}\otimes\bB{q}_{\alpha} \bigbcr_{i,j}
\bB{\eta}_{i,j}
- \frac{1}{2}  \theta^n_{0,i,j}\big( \bB{q}^n_{j,\alpha}
- \bB{q}^n_{i,\alpha}\big)  - \frac{1}{2}  \theta^n_{1,i,j}  \bB{\psi}^n_{i,j,\alpha}\\
& \qquad \qquad
- \frac{1}{2}  \theta^n_{1,i,j} \bigg(
\dfrac{1}{m^n_{j,\alpha}}\bB{q}^n_{j,\alpha}\otimes\bB{q}^n_{j,\alpha} -
\dfrac{1}{m^n_{i,\alpha}}\bB{q}^n_{i,\alpha}\otimes\bB{q}^n_{i,\alpha}
\bigg)\bB{\eta}_{i,j}
\,,
\end{align*}
where $b_{i,j,\alpha}^{n} := \bcl\rho_{\alpha}^n\bcr_{i,j}(z_{\mathrm{B},j}-z_{\mathrm{B},i})$
is included to preserve the well-balance and the stationary solution deduced in \cite[Proposition~1]{Osores2020}, and $\bB{\psi}_{i,j,\alpha}^n$ is approximated by averages and differences as follows
\begin{align*}
\begin{aligned}
 \bB{\psi}_{i,j,\alpha}^n & :=
 g\Biggl( \bigbcl m_{\alpha}^n\bigbcr_{i,j}
\Big(z_{\mathrm{B},j}-z_{\mathrm{B},i}+h_j^n-h_i^n\Big) +  \bigbcl(h^n)^2\bigbcr_{i,j} \biggl( \frac{l_{\alpha}}{2}
+\sum_{\beta=\alpha+1}^{M} l_{\beta} \biggr) \Big(\rho_{j,\alpha}^n - \rho_{i,\alpha}^n\Big)  \\
& \qquad \qquad  +\, \bigbcl h^n\bigbcr_{i,j}\sum_{\beta = \alpha+1}^{M} l_{\beta} \Big(m_ {j,\beta}^n-m_{i,\beta}^n-(m_{j,\alpha}^n-m_{i,\beta}^n)\Big)\Biggr)\bB{\eta}_{i,j}\,.
 \end{aligned}
\end{align*}
The coefficients $\theta^n_{0,i,j}$ and $\theta^n_{1,i,j}$ are respectively defined by
\begin{align*}
\theta^n_{0,i,j} = \frac{\sigma_{\mathrm{R},i,j}^n |\sigma_{\mathrm{L},i,j}^n|-\sigma_{\mathrm{L},i,j}^n|\sigma_{\mathrm{R},i,j}^n|}{
\sigma_{\mathrm{R},i,j}^n-\sigma_{\mathrm{L},i,j}^n}\,, \qquad
\theta^n_{1,i,j}= \frac{|\sigma_{\mathrm{R},i,j}^n|-|\sigma_{\mathrm{L},i,j}^n|}{\sigma_{\mathrm{R},i,j}^n-\sigma_{\mathrm{L},i,j}^n}\,,
\end{align*}
with $\smash{\sigma_{\mathrm{L},i,j}^n}$ and $\smash{\sigma_{\mathrm{R},i,j}^n}$ being the characteristic velocities, which are global approximations of the minimum and maximum wave speed of the viscosity matrix, respectively. We make use of the eigenvalues obtained in the one-dimensional case in \cite{Osores2019}, given by
\begin{align*}   %\label{wave1}
\sigma_{\mathrm{L},i,j}^n:=\bar{v}^n_{i,j}-(\chi\Psi_{i,j}^n)^{1/2},\qquad
\sigma_{\mathrm{R},i,j}^n:=\bar{v}^n_{i,j}+(\chi\Psi_{i,j}^n)^{1/2},
\end{align*}
where $\chi = 4-2/M$ is a positive constant and
\begin{align*}
\bar{v}^n_{i,j} & :=\frac{1}{M}\sum_{\beta=1}^M
\bigbcl\bB{\widetilde{v}}_{\beta}^{\:n}\bigbcr_{i,j}^{\tt t} \bB{\eta}_{i,j}\,, \qquad
\Psi_{i,j}^n  :=\sum_{\beta=1}^M(\bar{v}_{i,j}^n-v_{\beta,i,j}^n)^2 +
\frac{g}{2}\bigbcl h^n\bigbcr_{i,j} \biggl(1+\sum_{\beta=1}^M\dfrac{(2\beta-1)}{M\varrho_{\rm f}} \bigbcl\rho_{\beta}^n\bigbcr_{i,j} \biggr)\,,
\end{align*}
where $v_{\beta,i,j}^n := \bcl\bB{v}_{\beta}^n\bcr_{i,j}^{\tt t}\cdot \bB{\eta}_{i,j}$.
The horizontal fluxes corresponding to equations \eqref{eq:ml:c} and \eqref{eq:ml:s} are individually approximated using an upwind type scheme, these are respectively given by
\begin{align*}
\mathbcal{F}^{\bB{r},n}_{i,j,\alpha} &: = \upw\bigg(\mathcal{F}_{i,j,\alpha}^{m,n}\,;\, \dfrac{\bB{r}_{i,\alpha}^n}{m_{i,\alpha}^n},  \dfrac{\bB{r}_{j,\alpha}^n}{m_{j,\alpha}^n} \bigg),\qquad
\mathbcal{F}^{\bB{\zeta},n}_{i,j,\alpha}  := \upw\bigg(\mathcal{F}_{i,j,\alpha}^{m,n}\,;\, \dfrac{\bB{\zeta}_{i,\alpha}^n}{m_{i,\alpha}^n},  \dfrac{\bB{\zeta}_{j,\alpha}^n}{m_{j,\alpha}^n} \bigg)\,,
\end{align*}
where $\upw$ represents the upwind operator, which for a scalar quantity $\nu$, and vectors $\bB{\beta}$ and $\bB{\gamma}$ of the same length is defined by
\begin{align*}
\upw(\nu;\bB{\beta},\bB{\gamma}):= \nu \left( \frac{1+\sgn(\nu)}{2} \bB{\beta} +  \frac{1-\sgn(\nu)}{2} \bB{\gamma} \right)\,,
\end{align*}
with $\sgn$ represents the sign function.
In order to approximate the vertical fluxes of system \eqref{syst:finalmodel_ml}, which are determined by  \eqref{def:vfluxGi}, \eqref{def:vfluxSi} and \eqref{eq:fluxGnew}, we need to handle the convective and diffusive nonlinear fluxes $\bB{\tilde{f}}_{\alpha+1/2}$ and $\bB{\tilde{a}}_{\alpha+1/2}$ for $\alpha\in\intvl(M)$, respectively. To do so, we use the following approximations
\begin{align*}
\tilde{f}_{j,\alpha+1/2}^{(i),n}& := \frac{1}{2}
\Bigl(f_i\big(\bB{c}_{j,\alpha}^{n}\big) + f_i\big(\bB{c}_{j,\alpha+1}^{n}\big)\Bigr)
-\frac{\omega_{j,\alpha+1}}{2}\big(c_{j,\alpha+1}^{(i),n} - c_{j,\alpha}^{(i),n}\big)\\
&\quad -\frac{c_{j,\alpha}^{(i),n}}{2} \bigl|v_i^{\mathrm{MLB}}(\bB{c}_{j,\alpha+1}^n)-v_i^{\mathrm{MLB}}(\bB{c}_{j,\alpha}^n) \bigr|\sgn\big(c_{j,\alpha+1}^{(i),n}-c_{j,\alpha}^{(i),n}\big)\,,\quad \text{for } {i\in\intvl(n_{\bB{c}})}\,,\\
\bB{\tilde{a}}_{j,\alpha+1/2}^n &: = \dfrac{1}{2h_{\alpha}l_{\alpha}} \big(\mathbcal{D}(\bB{c}^n_{j,\alpha})+\mathbcal{D}(\bB{c}^n_{j,\alpha+1})\big)(\bB{c}_{j,\alpha+1}^{n}-\bB{c}_{j,\alpha}^{n})\,,
\end{align*}
where  $\omega_{\alpha+1} := \max_{i = 1,\dots,n_{\bB{c}}} |v_i^{\rm MLB}(\bB{c}_{\alpha+1})|$, and the matrix function $\mathbcal{D}$ is defined by its coefficients in \eqref{def:coeffD}. Then, for the case of zero flux boundary conditions at $z=h$ and $z = z_{\rm B}$, the fluxes in \eqref{def:vfluxGtilde} and \eqref{eq:fluxGnew} are respectively approximated by
\begin{align*}
\tilde{\mathcal{G}}_{i,j,\alpha+1/2}^n &
:=-\dfrac{\bcl\tilde{\rho}_{\alpha+1/2}^n\bcr_{i,j}}{\varrho_{\rm f}}\sum_{l=1}^{n_{\bB{c}}} \dfrac{(\varrho_l-\varrho_{\rm f})}{2\varrho_l}
\left(\tilde{f}_{i,\alpha+1/2}^{(l),n} + \tilde{f}_{j,\alpha+1/2}^{(l),n} - (\tilde{a}_{i,\alpha+1/2}^{(l),n}+\tilde{a}_{j,\alpha+1/2}^{(l),n})\right)\,,\\[1ex]
\mathcal{G}_{i,j,\alpha+1/2}^n &
= \dfrac{|e_{i,j}|}{|V_i|}\frac{\bcl\tilde{\rho}_{\alpha+1/2}^n\bcr_{i,j}}{\varrho_{\rm f}}  \sum_{\beta=1}^{\alpha} l_{\beta} \Big(  (\boldsymbol{R}_{j,\beta}^n-\bar{\boldsymbol{R}}_j^n)^{\tt t}-(\boldsymbol{R}_{i,\beta}^n-\bar{\boldsymbol{R}}_i^n)^{\tt t} \Big)\cdot \bB{\eta}_{i,j} + \tilde{\mathcal{G}}_{i,j,\alpha+1/2}^n\,.
\end{align*}
Then, the approximation of the vertical fluxes of equations \eqref{eq:ml:v}, \eqref{eq:ml:c} and \eqref{eq:ml:s} for all $\alpha\in \intvl_0(M)$, are respectively given by
\begin{align*}
\mathbcal{Q}_{i,\alpha+1/2}^n & := \sum_{j\in J_i} \bcl\bB{\tilde{q}}_{\alpha+1/2}^n\bcr_{i,j}\mathcal{G}_{i,j,\alpha+1/2}^n\,,\\[1ex]
%------------------------------------------------
\mathbcal{G}_{i,\alpha+1/2}^n & := \sum_{j\in J_i} \bcl\bB{\tilde{c}}_{\alpha+1/2}^n\bcr_{i,j}\mathcal{G}_{i,j,\alpha+1/2}^n -\Big(\bB{\tilde{f}}_{i,\alpha+1/2}^n-\bB{\tilde{a}}_{i,\alpha+1/2}^n\Big)\,,\\[1ex]
%%------------------------------------------------
% \begin{aligned}
\mathbcal{S}_{i,\alpha+1/2}^n & :=
\sum_{j\in J_i}\bcl\bB{\tilde{s}}_{\alpha+1/2}^n\bcr_{i,j}\mathcal{G}_{i,j,\alpha+1/2}^n\\
& \qquad \qquad+ {\rm Upw}\Bigg(\sum_{l=1}^{n_{\bB{c}}}\Big(\tilde{f}_{i,\alpha+1/2}^{(l),n}-\tilde{a}_{i,\alpha+1/2}^{(l),n}\Big); \dfrac{\bB{s}_{i,\alpha}^{n}}{\varrho_{\rm f}\phi_{{\rm f},i,\alpha}},
\dfrac{\bB{s}_{i,\alpha+1}^{n}}{\varrho_{\rm f}\phi_{{\rm f},i,\alpha+1}}\Bigg)\,.
% \end{aligned}\\
%%------------------------------------------------
 \end{align*}
 Finally, instead of approximating equation \eqref{eq:ml:m}, we determine the total mass $\bar{m}$ at each control volume from \eqref{eq:ml:mtot}, therefore the marching formula of the fully discrete finite volume scheme for all $V_i\in \mathcal{T}$, $\alpha\in\intvl(M)$ and $n\in \mathbb{N}$ reads as follows
 \begin{subequations} \label{syst:num}
\begin{align}
 \bar{m}_{i}^{n+1} & = \bar{m}_{i}^{n} -\dfrac{\Delta t}{|V_i|}\sum_{j\in J_i}\sum_{\beta=1}^{M} |e_{i,j}|l_{\beta}\mathcal{F}_{i,j,\beta}^{m,n} + \Delta t\sum_{\beta=1}^{M}l_{\beta}h_i^n \widetilde{\mathcal{R}}_{\rho,i,\beta}^n\,, \label{eq:num:mtot}\\[1ex]
 \bB{q}_{i,\alpha}^{n+1/2} & = \bB{q}_{i,\alpha}^n - \dfrac{\Delta t}{|V_i|}\sum_{j\in J_i} |e_{i,j}| \big(\mathbcal{F}^{\bB{q},n}_{i,j,\alpha}+\bB{\psi}_{i,j,\alpha}^n\big) + \dfrac{\Delta t}{l_{\alpha}}\big( \mathbcal{Q}_{i,\alpha+1/2}^n-\mathbcal{Q}_{i,\alpha-1/2}^n\big)\,,\label{eq:num:qhalf}\\[1ex]
 \bB{r}_{i,\alpha}^{n+1} & = \bB{r}_{i,\alpha}^n - \dfrac{\Delta t}{|V_i|}\sum_{j\in J_i} |e_{i,j}| \mathbcal{F}^{\bB{r},n}_{i,j,\alpha} + \dfrac{\Delta t}{l_{\alpha}}\big( \mathbcal{G}_{i,\alpha+1/2}^n-\mathbcal{G}_{i,\alpha-1/2}^n\big) + \Delta th_{i}^n\mathbcal{R}^{n}_{\bB{c},i,\alpha}\,,\label{eq:num:c}\\[1ex]
 \bB{\zeta}_{i,\alpha}^{n+1} & = \bB{\zeta}_{i,\alpha}^n -\dfrac{\Delta t}{|V_i|}\sum_{j\in J_i} |e_{i,j}| \mathbcal{F}^{\bB{\zeta},n}_{i,j,\alpha} + \dfrac{\Delta t}{l_{\alpha}}\big( \mathbcal{S}_{i,\alpha+1/2}^n-\mathbcal{S}_{i,\alpha-1/2}^n\big) + \Delta t h_{i}^n\mathbcal{R}^{n}_{\bB{s},i,\alpha}\,, \label{eq:num:s}
\end{align}
\end{subequations}
where $m_{i,\alpha}^{n+1}$ and $h_i^{n+1}$ are computed by
\begin{align}
 m_{i,\alpha}^{n+1} & = \bar{m}_{i}^{n+1}
 + \sum_{\beta=1}^{M}\sum_{l=1}^{n_{\bB{c}}} \dfrac{\varrho_l-\varrho_{\rm f}}{\varrho_l} l_{\beta}\big(r^{(l),n+1}_{i,\alpha}-r^{(l),n+1}_{i,\beta}\big)\,, \label{eq:update:m}\\[1ex]
 h_i^{n+1} & = \dfrac{1}{\varrho_{\rm f}}\biggl(\bar{m}_{i}^{n+1}
 - \sum_{\beta=1}^{M}\sum_{l=1}^{n_{\bB{c}}} \dfrac{\varrho_l-\varrho_{\rm f}}{\varrho_l} l_{\beta}r^{(l),n+1}_{i,\beta}
 \biggr)\,. \label{eq:update:h}
\end{align}
To determine the velocity field at time $t^{n+1}$, we need to include the discretization of the viscous terms \eqref{def:Kflux}, this is done by solving the following linearly implicit equation
\begin{align}\label{eq:update:q}
 \bB{q}_{i,\alpha}^{n+1} =
 \bB{q}_{i,\alpha}^{n+1/2}
 + \dfrac{\Delta t}{l_{\alpha}}\Big(\bB{K}_{i,\alpha+1/2}^{n+1}  - \bB{K}_{i,\alpha-1/2}^{n+1} \Big)\qquad \forall \alpha\in \intvl(M)\,,
\end{align}
where
\begin{align*}
\bB{K}_{i,\alpha+1/2}^{n+1} := \dfrac{\mu(\phi_{i,\alpha+1/2}^{n+1})}{2(\tfrac{1}{2}l_{\alpha+1} + \tfrac{1}{2}l_{\alpha}) h_i^n}\bigg(\dfrac{\bB{q}_{i,\alpha+1}^{n+1}}{m_{i,\alpha+1}^{n+1}} - \dfrac{\bB{q}_{i,\alpha}^{n+1}}{m_{i,\alpha}^{n+1}}\bigg)\qquad \forall\alpha\in \intvl_0(M)\,,
\end{align*}
with $\phi_{i,\alpha+1/2}^n:= \tfrac{1}{2}\big(\phi(\bB{c}_{i,\alpha+1}^n) + \phi(\bB{c}_{i,\alpha}^n)\big)$. Note that equation~\eqref{eq:update:q} can be written as a linear system composed by an invertible tri-diagonal matrix. The numerical scheme is supplemented with the $\mathrm{CFL}$ condition provided in \cite{Osores2020}, which corresponds to an adaptive time-stepping given by
\begin{align*}%\label{eq:cfl}
\max\left\{\frac{|\lambda_{i,j}|}{d_{i,j}} : \quad i\in \intvl(N)\,,\quad j\in J_i\,\right\} \Delta t\,\leq\, \mathrm{CFL}\,,
\end{align*}
where $\lambda_{i,j}$ are bounds of the eigenvalues of the viscosity matrix,
$d_{i,j}=||\bB{x}_j-\bB{x}_i||_2$ is the distance between centers of the volume $i$ and $j$. For all examples, we have set the Courant number as $\mathrm{CFL}=0.5$.

\section{Numerical simulations}\label{sec:numexamples}

We have implemented the numerical scheme composed by \eqref{syst:num}, \eqref{eq:update:m},  \eqref{eq:update:h} and \eqref{eq:update:q} in the programming language Fortran 90, and the tri-diagonal linear system arising from \eqref{eq:update:q} is computed using the library SuperLU \cite{Li05}.
For the numerical examples, we consider the reduced denitrification process \cite{SDcace_reactive}, which consists of converting nitrate into nitrogen gas. The solid phase is composed by two species, the ordinary heterotrophic organisms $c^{(1)}$ and undegradable organics $c^{(2)}$, and the substrates are three, the diluted nitrate $s^{(1)}$, the readily biodegradable substrate $s^{(2)}$, and the diluted nitrogen $s^{(3)}$. The vectors of concentration  are $\bB{c} = \big(c^{(1)},c^{(2)}\big)$ and  $\bB{s} = \big(s^{(1)},s^{(2)},s^{(3)}\big)$, which in \cite{SDcace_reactive} are denoted by $\bB{c}=\big(X_{\rm OHO},X_{\rm U}\big)$ and $\bB{s}=\big(S_{\rm NO_3},S_{\rm S}, S_{\rm N_2}\big)$. The vector of reaction rates and stoichiometric matrices are respectively given by
\begin{align*}
 \bB{\kappa}(\bB{c},\bB{s}) = c^{(1)}\begin{bmatrix}
           \mu(s^{(1)},s^{(2)})\\
           b
          \end{bmatrix}\,,\qquad
  \boldsymbol{\sigma}_{\boldsymbol{c}} &=
\begin{bmatrix}
 1 & -1\\
 0 & f_{\rm p}
\end{bmatrix},\qquad \bB{\sigma}_{\bB{s}} = \begin{bmatrix}
-\bar{Y} & -1/Y & \bar{Y}\\
 0 & 0.8 & 0
\end{bmatrix}\,,
\end{align*}
where $b=6.94\times10^{-5}{\rm s^{-1}}$ is the decay rate of heterotrophic organisms, $f_{\rm p} = 0.2$ is the portion of these that decays to undegradable organics, $Y=0.67$ and $\bar{Y} = 0.172216$ are (dimensionless) yield factors.  The specific growth rate function is defined by
\begin{align*}
 \mu(\bB{s}) = \mu_{\max} \dfrac{s^{(1)}}{\kappa_1+s^{(1)}}\dfrac{s^{(2)}}{\kappa_2+s^{(2)}}\,,
\end{align*}
with $\mu_{\max} = 5.56\times 10^{-4}\,{\rm s}^{-1}$ being the maximum growth rate, and $\kappa_1=5\times 10^{-4}\,{\rm kg/m^3}$ and $\kappa_2=0.02\,{\rm kg/m^3}$ are saturation constants. Note that constants $\mu_{\max}$ and $b$ are chosen $10$ times larger than those used in the reduced denitrification model \cite{SDcace_reactive,SDIMA_MOL}, this is in order to amplify the effect produced by the reaction terms in a shorter simulation time. Then, the reaction terms read as follows
\begin{align*}
 \mathbcal{R}_{\bB{c}}(\bB{c},\bB{s}) &=
 c^{(1)}\begin{bmatrix}
                                          \mu(\bB{s}) - b\\[0.5ex]
                                          f_{\rm p}b
                                         \end{bmatrix}\,,\qquad
 \mathbcal{R}_{\bB{s}}(\bB{c},\bB{s}) = c^{(1)}\begin{bmatrix}
                                          -\bar{Y}\mu(\bB{s})\\[0.5ex]
                                          -\mu(\bB{s})/Y + 0.8b\\[0.5ex]
                                          \bar{Y}\mu(\bB{s})
                                         \end{bmatrix}\,,
\end{align*}
The hindered settling, effective solid stress and viscosity functions used in all simulations are provided respectively as follows
\begin{align*}
 v_{\rm hs}(\phi) := \bigg(1 - \dfrac{\phi}{\phi_{\max}}\bigg)^{4.7}\,,\qquad
\sigma_{\rm e}(\phi) :=
\begin{cases}
\sigma_0 \big((\phi / \phi_{\rm c})^{5.0} - 1\big) & \mbox{if }\phi \geq \phi_{\rm c}\,,\\
0 & \mbox{if } \phi < \phi_{\rm c}\,,
\end{cases}
\end{align*}
and
\begin{align*}
\mu(\phi) :=  \mu_0\bigg(1 - \dfrac{\phi}{0.95}\bigg)^{-2.5}\,,
\end{align*}
where $\sigma_0 = 0.02\, \rm m^2/s^2$, $\mu_0=0.01\,\rm Pa\,s$, and the maximal and critical volume fractions are set to $\phi_{\max} = 0.02$ and $\phi_{\rm c}:=0.003$, respectively. Other parameters are $\varrho_1 = \varrho_2 = 2000\,\rm kg/m^3$, $\varrho_{\rm f} = 998\,\rm kg/m^3$, $\delta_1 = 4\times 10^{-4}$, $\delta_2 = 2.5\times 10^{-4}$ and $g = 9.81\,\rm m/s^2$. In all simulations we use the following homogeneous initial conditions for the substrates
\begin{align*}
 \bB{s}(\bB{x},0) = (0.006, 0.0009, 0.0)^{\tt t}\,{\rm kg/m^3}\,\qquad\forall \bB{x}\in\Omega\,,
\end{align*}
and for all simulations except the second one, we use the solids initial condition
\begin{align*}
 \bB{c}(\bB{x},0) = (3,2.5)^{\tt t}\,{\rm kg/m^3}\,\qquad\forall \bB{x}\in\Omega\,.
\end{align*}
Finally, we implemented zero flux-boundary conditions so that $\widetilde{\bB{v}}_\alpha^{\tt t}\cdot\bB{\eta}_{\alpha} = 0$ on $\partial\Omega(t)$, for all $\alpha\in\intvl(M)$, where $\bB{\eta}_{\alpha}$ is the outward-facing normal to the boundary $\partial\Omega_{\alpha}$.
% \pagecolor{red}
\newcommand{\escala}{0.55}
\begin{figure}[t]
 \begin{tabular}{cc}
 \multicolumn{2}{c}{\bf Ordinary heterotrophic organisms ($X_{\rm OHO}$)}\\[1.5ex]
 \includegraphics[scale=\escala]{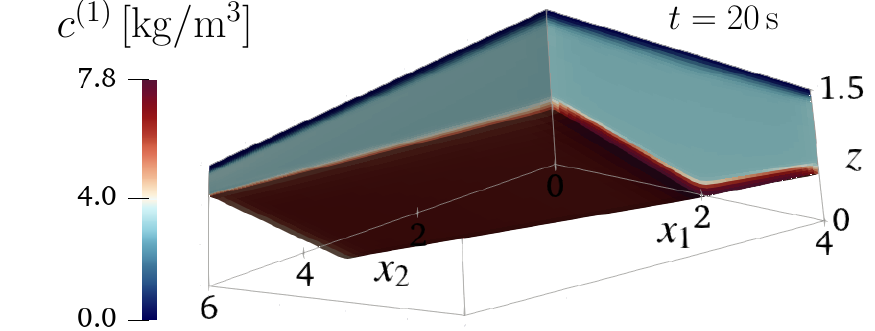} & \includegraphics[scale=\escala]{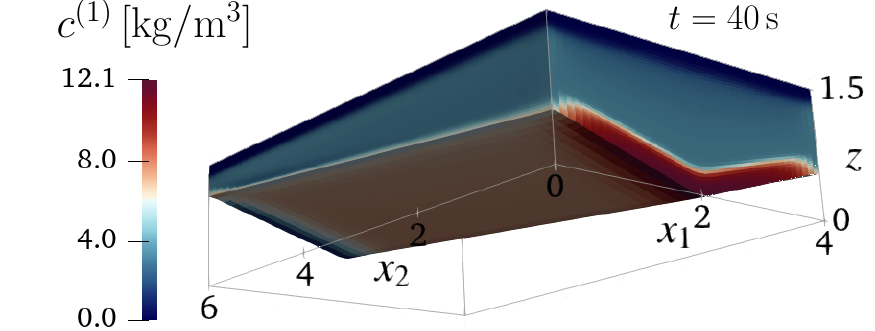}\\[2ex]
 \multicolumn{2}{c}{\bf Undegradable organics ($X_{\rm U}$)}\\[1.5ex]
 \includegraphics[scale=\escala]{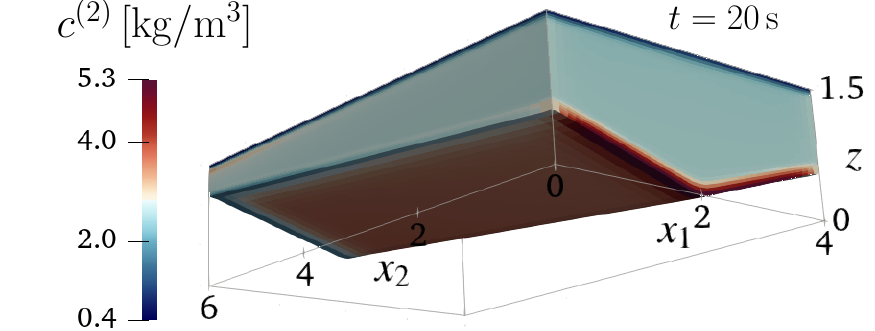} & \includegraphics[scale=\escala]{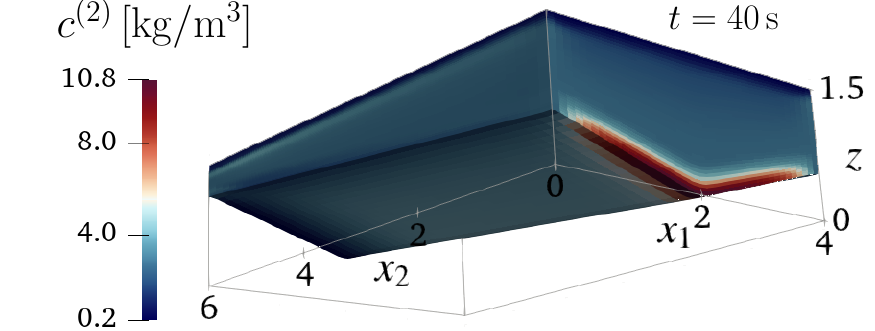}\\[2ex]
  \multicolumn{2}{c}{\bf Mass-average velocity}\\[1.5ex]
  \includegraphics[scale=\escala]{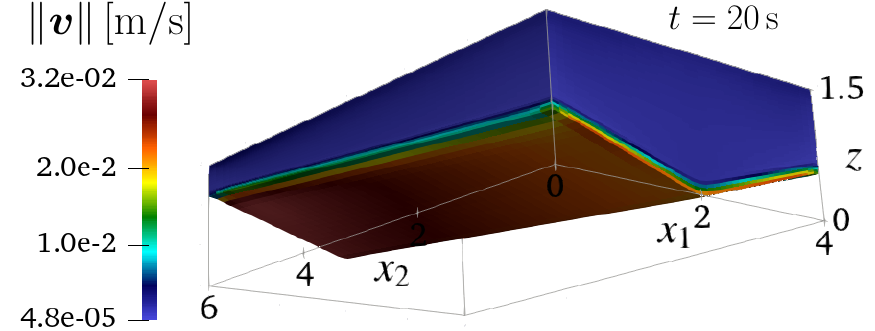} & \includegraphics[scale=\escala]{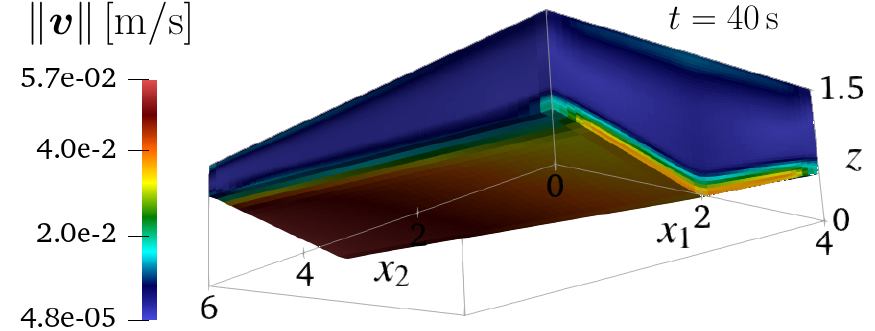}\\
 \end{tabular}
 \caption{Simulation 1. Concentration of solid components and magnitude of the mass-average velocity at a middle time $t=20\,\rm s$ (first column), and final time $t=40\,\rm s$ (second column).}\label{fig:sim1:1}
\end{figure}

\begin{figure}[t]
 \begin{tabular}{cc}
 \multicolumn{2}{c}{\bf Nitrate substrate ($S_{\rm NO_3}$)}\\[1.5ex]
 \includegraphics[scale=\escala]{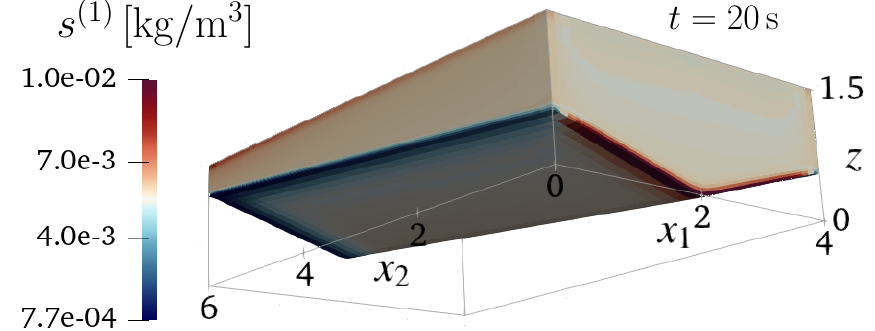} & \includegraphics[scale=\escala]{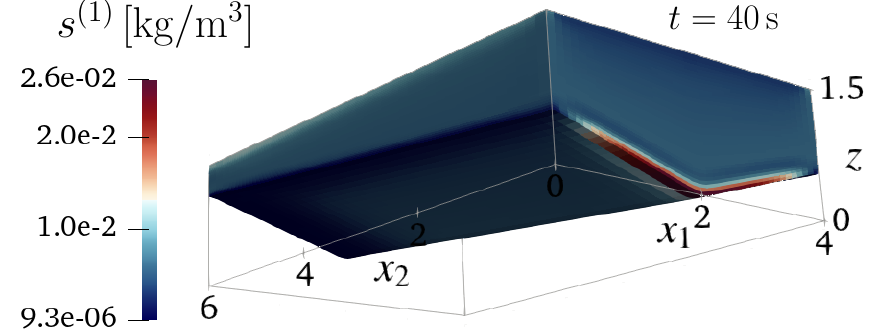}\\[2ex]
  \multicolumn{2}{c}{\bf Readily biodegradable substrate ($S_{\rm S}$)}\\[1.5ex]
 \includegraphics[scale=\escala]{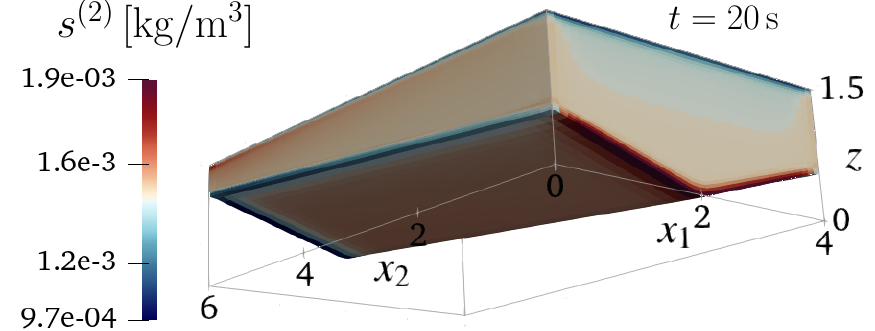} & \includegraphics[scale=\escala]{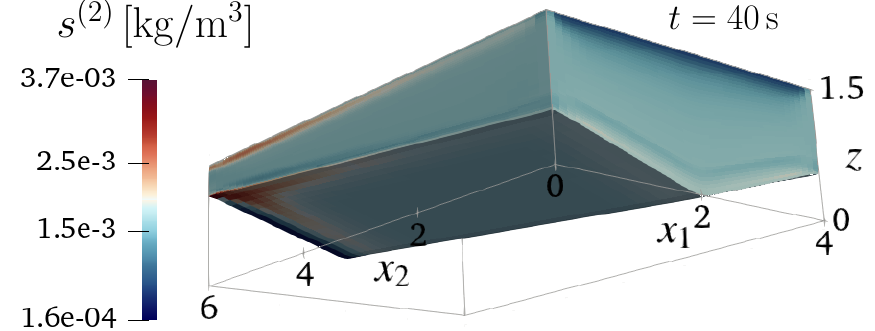}\\[2ex]
   \multicolumn{2}{c}{\bf Nitrogen substrate ($S_{\rm N_2}$)}\\[1.5ex]
 \includegraphics[scale=\escala]{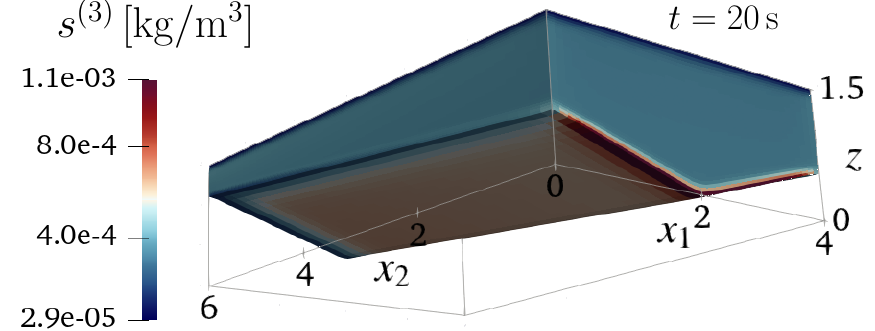} & \includegraphics[scale=\escala]{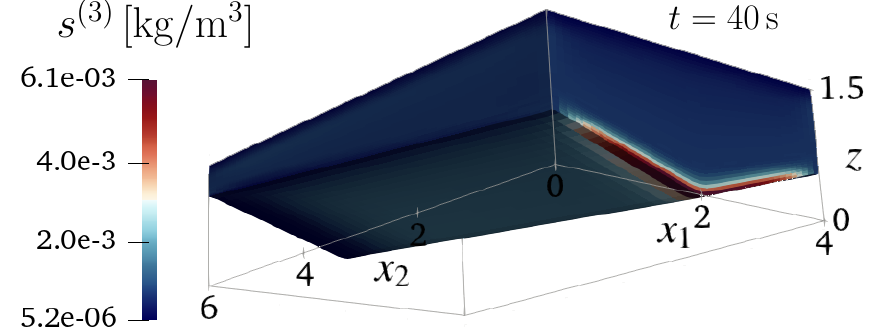}
 \end{tabular}
 \caption{Simulation 1. Concentration of substrates at a middle time $t=20\,\rm s$ (first column), and final time $t=40\,\rm s$ (second column).}\label{fig:sim1:2}
\end{figure}
\subsection{Simulation 1: Inclined walls settling}
\begin{figure}[t]
 \begin{tabular}{cc}
  \multicolumn{2}{c}{\bf Ordinary heterotrophic organisms ($X_{\rm OHO}$)}\\[1.5ex]
 \includegraphics[scale=\escala]{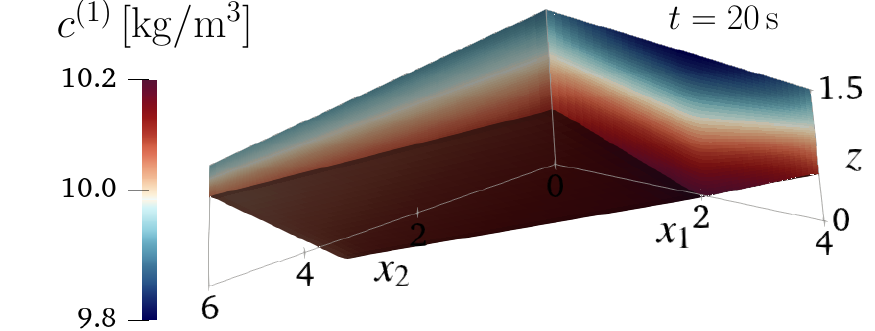} & \includegraphics[scale=\escala]{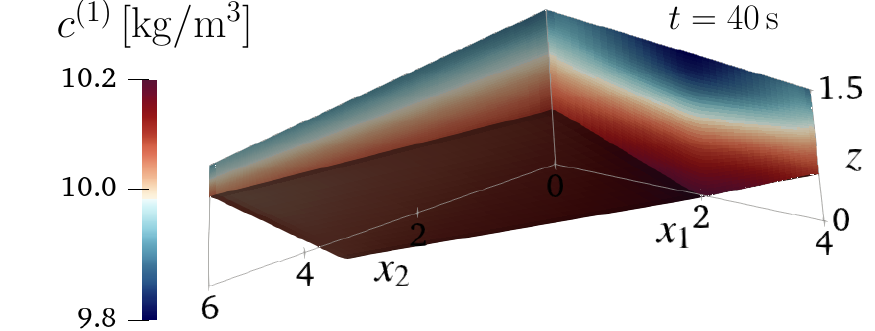}\\[2ex]
  \multicolumn{2}{c}{\bf Undegradable organics ($X_{\rm U}$)}\\[1.5ex]
 \includegraphics[scale=\escala]{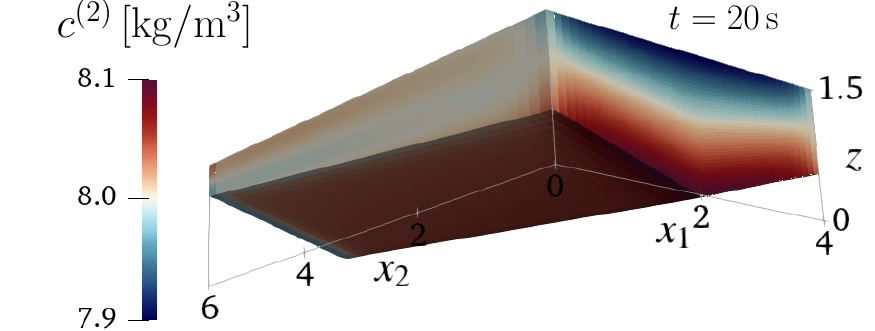} & \includegraphics[scale=\escala]{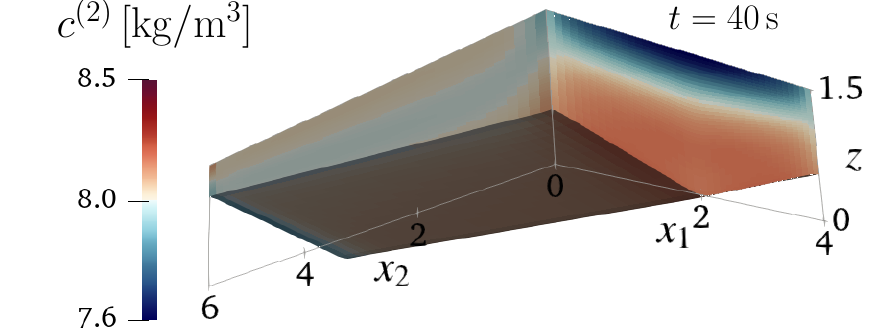}\\[2ex]
   \multicolumn{2}{c}{\bf Mass-average velocity}\\[1.5ex]
  \includegraphics[scale=\escala]{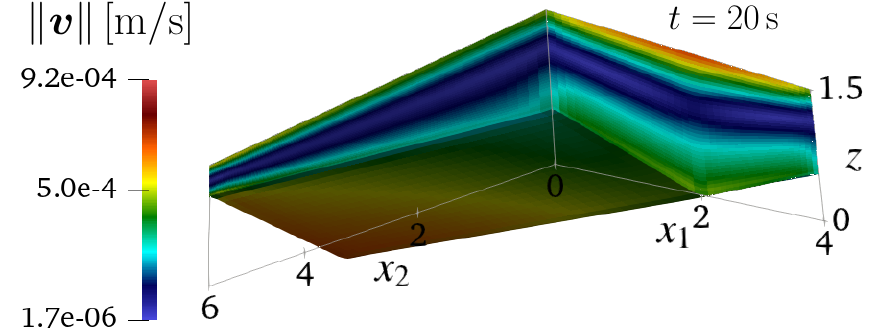} & \includegraphics[scale=\escala]{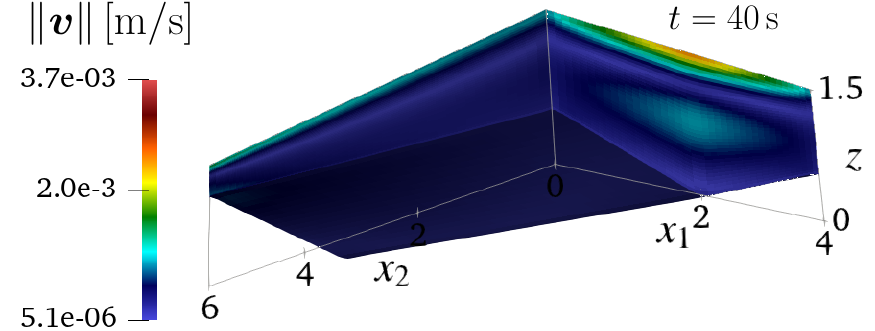}\\
 \end{tabular}
 \caption{Simulation 2. Concentration of solid components and magnitud of the mass-average velocity at a middle time $t=20\,\rm s$ (first column), and final time $t=40\,\rm s$ (second column).}\label{fig:sim2:1}
\end{figure}

\begin{figure}[t]
 \begin{tabular}{cc}
 \multicolumn{2}{c}{\bf Nitrate substrate ($S_{\rm NO_3}$)}\\[1.5ex]
 \includegraphics[scale=\escala]{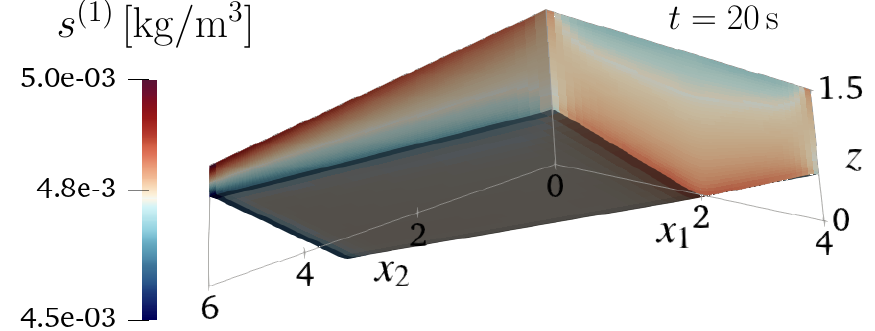} & \includegraphics[scale=\escala]{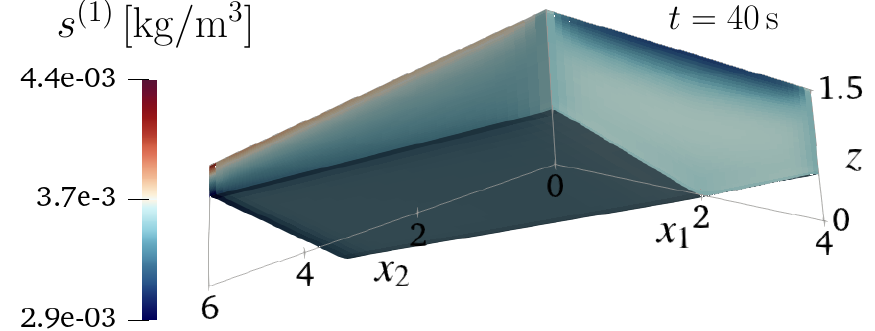}\\[2ex]
\multicolumn{2}{c}{\bf Readily biodegradable substrate ($S_{\rm S}$)}\\[1.5ex]
 \includegraphics[scale=\escala]{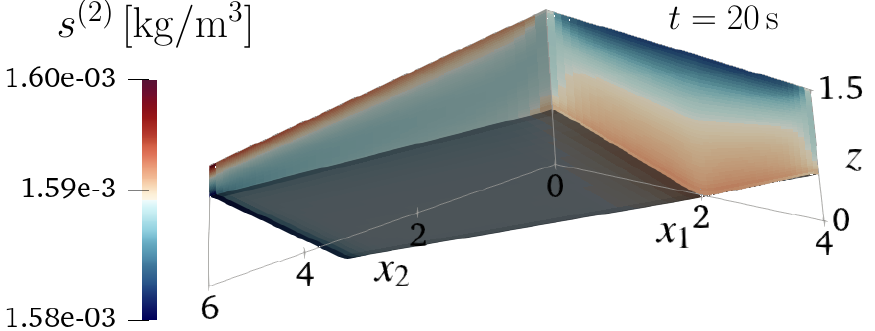} & \includegraphics[scale=\escala]{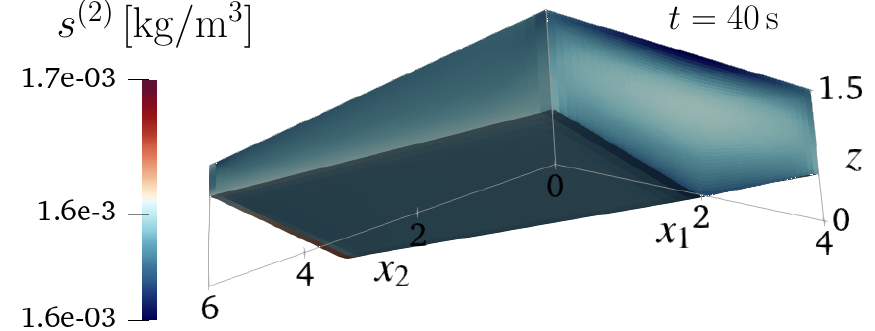}\\[2ex]
\multicolumn{2}{c}{\bf Nitrogen substrate ($S_{\rm N_2}$)}\\[1.5ex]
 \includegraphics[scale=\escala]{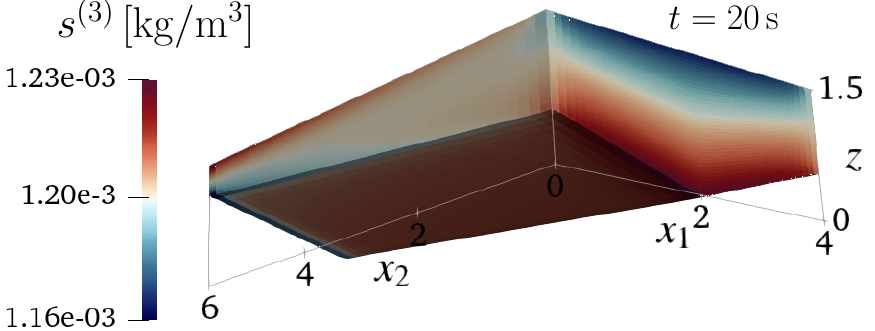} & \includegraphics[scale=\escala]{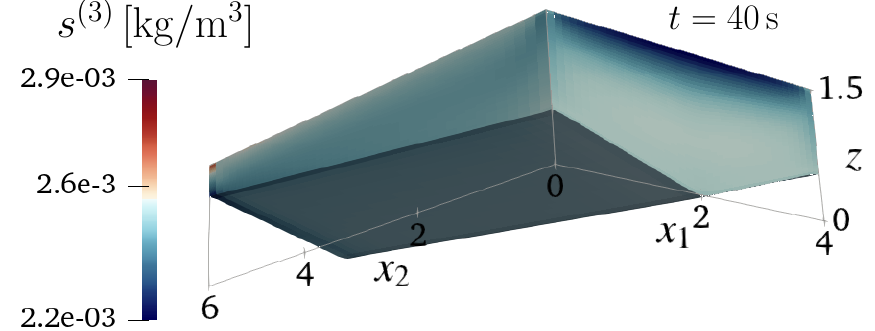}
 \end{tabular}
 \caption{Simulation 2. Concentration of substrates at a middle time $t=20\,\rm s$ (first column), and final time $t=40\,\rm s$ (second column).} \label{fig:sim2:2}
\end{figure}

We consider a closed vessel with an inclined bottom wall, enclosed on a rectangle of $4\,\rm m$ wide ($x_1$ direction) and $6\,\rm m$ long ($x_2$ direction), such that the projected domain corresponds to $\bB{\Pi}(0) = (0,4)\times (0,6)$.
The vessel, shown in Figures~\ref{fig:sim1:1} to \ref{fig:sim2:2}, has a bottom wall described by
\begin{align*}
 z_{\rm B}(\widetilde{\bB{x}}) = \max\bigg\{\pm\dfrac{6.4125}{22.8}(x_1 - 2) + \dfrac{2.25}{22.8}(x_2 - 6) + 0,5625, 0\bigg\}\,,\qquad \widetilde{\bB{x}}\in \bB{\Pi}(0)\,,
\end{align*}
which is initially filled with mixture up to $z = 1.5\,\rm m$, the initial height is $h(\widetilde{\bB{x}},0) = 1.5\,\rm m$, for $\widetilde{\bB{x}}\in \bB{\Pi}(0)$.
We set the mesh sizes $\Delta x = 0.1\,\rm m$ and $\Delta y = 0.15\,\rm m$, and $M=40$ layers, which results in a total of $64000$ cells, counting all layers. Figure~\ref{fig:sim1:1} shows the simulated concentration of solid species $c^{(1)}$ and $c^{(2)}$, and norm of the mass-average velocity at times $t=20\,\rm s$ and $t=T=40\,\rm s$ (the final simulation time). Both components $c^{(1)}$ and $c^{(2)}$ settle down on the bottom wall, and due to the inclination of the vessel, the sedimentation continues towards $x_2=0$ reaching maximum values of concentration near the center line $x_1=2$. The latter proofs that the model and numerical scheme properly address the movements of the particles in all three directions.
Moreover, the mass-average velocity, shown in the third row of Figure~\ref{fig:sim1:1}, assumes its maximum at the bottom layers, remaining almost constant in the rest of the domain. The above is explained by the fact that a greater amount of solid particles are being concentrated in the first layers, inducing the movement along them. The concentration of substrate components $s^{(1)}$, $s^{(2)}$ and $s^{(3)}$ at both times are presented in Figure~\ref{fig:sim1:2}. For the nitrate substrate $s^{(1)}$ (first row, Figure~\ref{fig:sim1:2}), we observe that the maximum concentration is assumed at regions where $c^{(1)}$ is higher. However, at time $t=20\,\rm s$ a thin layer of this substrate is still remaining at the top of the vessel. On the other hand, the readily biodegradable substrate $s^{(2)}$ (second row of Figure~\ref{fig:sim1:2}), at time $t=20\,\rm s$ has increased in the entire domain, and at time $t=40\,\rm s$ exhibit a decreasing behavior near the top of the wall at $x_2=0$. Furthermore, the maximum concentration is located near the corners at $x_2=6$. The nitrogen substrate $s^{(3)}$ (third row, Figure~\ref{fig:sim1:2}), which is entirely produced by the chemical reactions, is accumulated with greater abundance at the bottom wall towards $x_2=0$, and it increases over time.

\subsection{Simulation 2: Higher compression effect}
We consider the same set of parameters and domain as in Simulation 1, but amplify the compression effect by setting $\sigma_0 =  0.5\,\rm m^2/s^2$ and $\phi_{\rm c}= 0.002$. We also use a higher initial condition for the solids
\begin{align*}
 \bB{c}(\bB{x},0) = (10,8)^{\tt t}\,{\rm kg/m^3}\qquad \forall \bB{x}\in \Omega\,.
\end{align*}
In Figure~\ref{fig:sim2:1}, we report the simulated concentration of the solid components $c^{(1)}$ and $c^{(2)}$ at times $t=20\,\rm s$ and $t=40\,\rm s$. As expected, sedimentation evolves at a slower rate, with smoother space variations within the vessel. In both components, solid particles are accumulated downward and towards the wall at $x_2 = 0$ as in Simulation 1, which shows that higher compression does not prevent horizontal movements. The velocity field (third row, Figure~\ref{fig:sim2:1}) differs considerably to the one in Simulation~1. At $t=20\,\rm s$, two regions of higher velocity appear at the top and bottom of the vessel, and then at $t=40\,\rm s$, the velocity is maximal in the top layer. The later is explained since the vessel is full of highly concentrated solid particles and no separation from clear water has occurred yet. Compared to Simulation~1, the three substrates show an increase in their concentrations (Figure~\ref{fig:sim2:2}), where nitrate $s^{(1)}$ presents the greatest variation. The three substrates assume their maximum values near the top of the wall at $x_2=6$.

\newcommand{\escalan}{0.6}
\begin{figure}[!t]
\centering
 \begin{tabular}{cc}
 {\bf Solid $X_{\rm OHO}$} & {\bf Substrate $S_{\rm NO_3}$}\\
 \includegraphics[scale=\escalan]{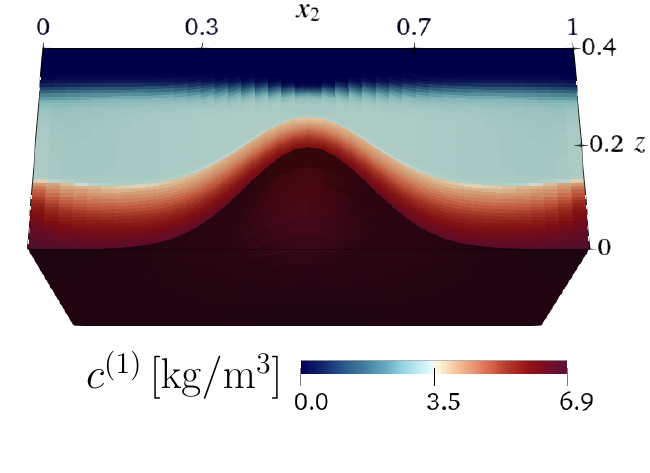} & \includegraphics[scale=\escalan]{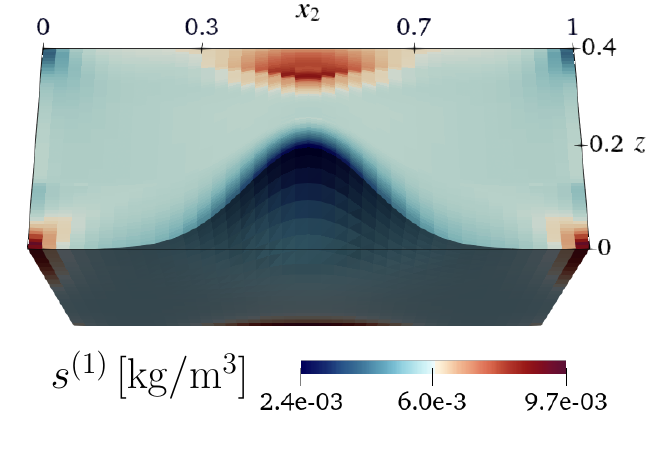}\\
  {\bf Solid $X_{\rm U}$} & {\bf Substrate $S_{\rm S}$}\\
 \includegraphics[scale=\escalan]{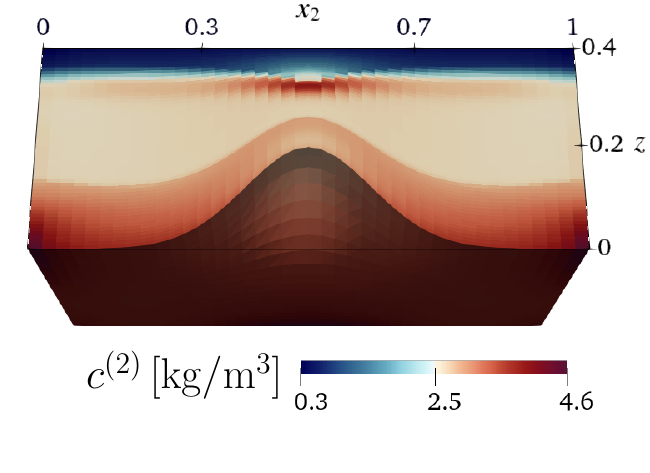} & \includegraphics[scale=\escalan]{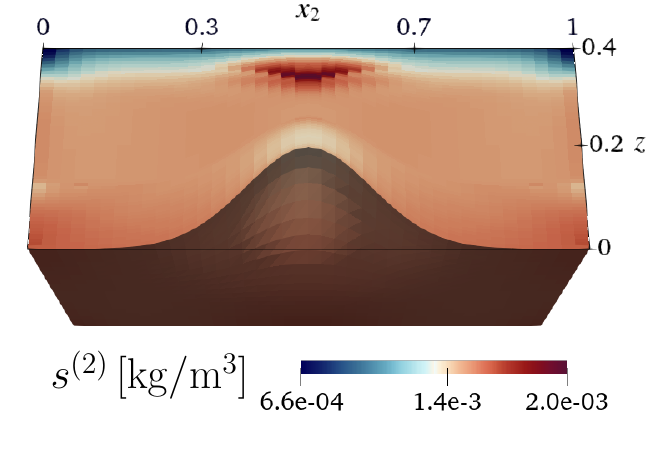}\\
  {\bf Velocity} & {\bf Substrate $S_{N_2}$}\\
  \includegraphics[scale=\escalan]{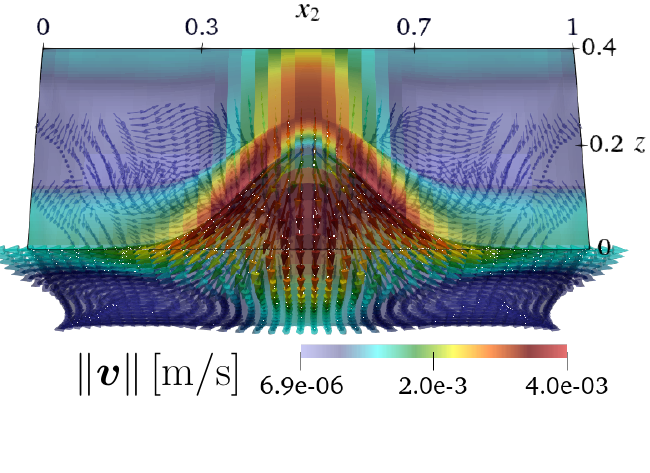} & \includegraphics[scale=\escalan]{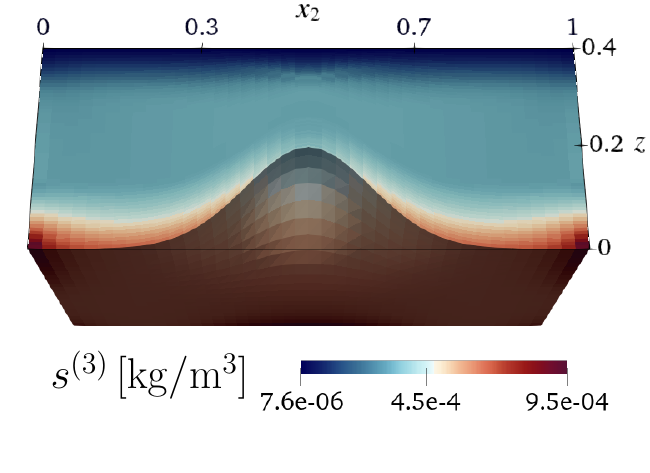}\\
 \end{tabular}
 \caption{Simulation 3. Concentration of the solid and substrate components, and velocity field at the end time $t=20\,\rm s$ on the sub-domain for $0\leq x_1\leq 0.5$. The plot of the velocity field includes the arrows of its direction.}\label{fig:sim3}
\end{figure}
\subsection{Simulation 3: Concave bathymetry}

We assume now that the sedimentation process is carried out in a domain with concave and smooth bathymetry. For this end, we consider the square projected domain $\bB{\Pi}(0) = (0,1)^2$ and the following bathymetry
\begin{align*}
 z_{\rm B}(\widetilde{\bB{x}}) = 0.2\exp\Big(-40\big((x_1-0.5)^2 + (x_2-0.5)^2\big)\Big)\qquad \forall \widetilde{\bB{x}}\in \bB{\Pi}(0)\,.
\end{align*}
For this simulation we set $\Delta x = \Delta y = 0.025\,\rm m$, final simulation time $T=20\,\rm s$, total number of layers $M = 60\,\rm s$, and initial height $h(\widetilde{\bB{x}},0) = 0.4\,\rm m$.
In Figure~\ref{fig:sim3}, we present the numerical results of the concentrations of each solid  species and substrates, and velocity field in a domain that is cut off at $x_1=0.5$ (for $0\leq x_1\leq 0.5$).
We observe that the concentration of particles $c^{(1)}$ accumulates at the bottom and towards the side walls, while particles $c^{(2)}$ remain floating at the interface between the solid and liquid phase. The velocity field assumes its maximum near the top of the curved part of $z_{\rm B}$ and in the region of lower height $h$, and recirculation of the fluid is observed. Concentration $s^{(1)}$ accumulates at the top near $\widetilde{\bB{x}} = (0,0)$ and at the bottom-central part of the side walls. Component $s^{(2)}$ is less concentrated at the top of the domain, with the exception of a central region of highly concentrated substrate near $\widetilde{\bB{x}} = (0,0)$. Substrate $s^{(3)}$ exhibits behavior more similar to that of $c^{(1)}$ with a thinner layer of concentrated substrate and maximum towards the bottom of the side walls.

\section{Conclusions}\label{sec:conclusions}

We have extended the one-dimensional reactive sedimentation model provided in \cite{SDIMA_MOL} to the three-dimensional case combining the multilayer shallow water approach for polydisperse settling presented in \cite{Osores2019,Osores2020}. We consider solid particles featuring varying diameters and densities which settle with different velocities, while the substrates or liquid components move with the same velocity.
The governing equations are written in terms of concentrations of solid species and substrates, quantities that are interconnected through nonlinear reaction terms, as well as by their respective velocity fields.
Moreover, the model includes nonlinear convective velocities involving the hindered settling function, and the compressibility of the sediment, both acting on the vertical component. For the multilayer version of the model, the momentum equation which is written in terms of the horizontal components of the mass-average velocity, we have assumed hydrostatic pressure and a reduced version of the viscous stress tensor. One of the advantages of introducing this multilayer approach, is that at each layer, we only have horizontal space derivatives, while vertical fluxes and reaction terms are considered as source terms. In addition, the total height of the fluid column allows determining the evolution of the free surface.

Regarding the numerical scheme, we approximate the vertical fluxes in a slightly different way than in \cite{Osores2020}, where the approximation of the vertical fluxes corresponding to the substrates depend on the numerical fluxes of each solid species. Numerical examples show the good performance of the model and numerical scheme in scenarios varying the compression effect and bathymetry for the denitrification model~\cite{SDcace1}. We emphasize that our work can be used to simulate more general processes involving concentrations of solid particles and substrates, such as the ASM1 \cite{Henze2000ASMbook} or eutrophication of shallow lakes \cite{Jayaweera1996, Jorgensen1988, Pauer2000, vanderMolen1994}.

Finally, there remain unresolved mathematical questions concerning the model equations, including the hyperbolicity of the system and the existence of entropy solutions. Additionally, questions persist regarding the analysis of the numerical scheme, encompassing aspects such as positivity preserving, stability, and convergence analysis. Further studies can be done in order to determine the appropriateness of neglecting the vertical flux contributions coming from the reaction terms, and in the line of considering alternative  numerical schemes such as the discontinuous Galerkin and entropy stable method proposed in \cite{Wintermeyer2017}.

\section*{Acknowledgments} JC is supported by ANID-Chile through Fondecyt Postdoctoral project No.~3230553.

\bibliographystyle{plain}
\bibliography{ref_shw}

\begin{thebibliography}{10}

\bibitem{Anestis1981}
G.~Anestis.
\newblock {\em Eine eindimensionale {T}heorie der {S}edimentation in
  {A}b\-setz\-beh\"{a}ltern ver\"ander\-lichen {Q}uerschnitts und in
  {Z}entrifugen}.
\newblock PhD thesis, TU Vienna, Austria, 1981.

\bibitem{Audusse2005}
E.~Audusse.
\newblock A multilayer saint-venant model: Derivation and numerical validation.
\newblock {\em Discrete Continuous Dyn.~Syst.~Ser.~B}, 5(2):189--214, 2005.

\bibitem{Audusse2011}
E.~Audusse, M.-O. Bristeau, M.~Pelanti, and J.~Sainte-Marie.
\newblock Approximation of the hydrostatic navier{\textendash}stokes system for
  density stratified flows by a multilayer model: Kinetic interpretation and
  numerical solution.
\newblock {\em J.~Comput.~Phys.}, 230(9):3453--3478, 2011.

\bibitem{Audusse2010}
E.~Audusse, M.-O. Bristeau, B.~Perthame, and J.~Sainte-Marie.
\newblock A multilayer saint-venant system with mass exchanges for shallow
  water flows. derivation and numerical validation.
\newblock {\em ESAIM: Math.~Model.~Numer.~Anal.}, 45(1):169--200, 2010.

\bibitem{Basson2009}
D.K. Basson, S.~Berres, and R.~B\"{u}rger.
\newblock On models of polydisperse sedimentation with particle-size-specific
  hindered-settling factors.
\newblock {\em Appl. Math. Model.}, 33(4):1815--1835, 2009.

\bibitem{SDsiap_varyingA}
R.~B\"{u}rger, J.~Careaga, and S.~Diehl.
\newblock Entropy solutions of a scalar conservation law modeling sedimentation
  in vessels with varying cross-sectional area.
\newblock {\em {SIAM} J. Appl. Math.}, 77(2):789--811, 2017.

\bibitem{SDcec_varyingA}
R.~B\"urger, J.~Careaga, and S.~Diehl.
\newblock A simulation model for settling tanks with varying cross-sectional
  area.
\newblock {\em Chem. Eng. Commun.}, 204(11):1270--1281, 2017.

\bibitem{SDIMA_MOL}
R.~B\"urger, J.~Careaga, and S.~Diehl.
\newblock A method-of-lines formulation for a model of reactive settling in
  tanks with varying cross-sectional area.
\newblock {\em IMA J. Appl. Math.}, 86(3):514--546, 2021.

\bibitem{SDcace_reactive}
R.~B{\"u}rger, J.~Careaga, S.~Diehl, C.~Mej\'{i}as, I.~Nopens, E.~Torfs, and
  P.~A. Vanrolleghem.
\newblock Simulations of reactive settling of activated sludge with a reduced
  biokinetic model.
\newblock {\em Computers Chem. Eng.}, 92:216--229, 2016.

\bibitem{SDAMM_SBR1}
R.~B\"{u}rger, J.~Careaga, S.~Diehl, and R.~Pineda.
\newblock A moving-boundary model of reactive settling in wastewater treatment.
  {P}art~1: Governing equations.
\newblock {\em Appl. Math. Modelling}, 106:390--401, 2022.

\bibitem{SDAMM_SBR2}
R.~B\"{u}rger, J.~Careaga, S.~Diehl, and R.~Pineda.
\newblock A moving-boundary model of reactive settling in wastewater treatment.
  {P}art~2: {N}umerical scheme.
\newblock {\em Appl. Math. Modelling}, 111:247--269, 2022.

\bibitem{Burger2023a}
R.~B\"{u}rger, J.~Careaga, S.~Diehl, and R.~Pineda.
\newblock A model of reactive settling of activated sludge: Comparison with
  experimental data.
\newblock {\em Chem.~Eng.~Sci.}, 267:118244, March 2023.

\bibitem{Burger2023b}
R.~B\"{u}rger, J.~Careaga, S.~Diehl, and R.~Pineda.
\newblock Numerical schemes for a moving-boundary convection-diffusion-reaction
  model of sequencing batch reactors.
\newblock {\em ESAIM: Math.~Model.~Numer.~Anal.}, 57(5):2931--2976, 2023.

\bibitem{Burger&D&K2004}
R.~B{\"u}rger, J.~J.~R. Damasceno, and K.~H. Karlsen.
\newblock A mathematical model for batch and continuous thickening of
  flocculated suspensions in vessels with varying cross-section.
\newblock {\em Int. J. Miner. Process.}, 73:183--208, 2004.

\bibitem{SDcace1}
R.~B{\"u}rger, S.~Diehl, S.~Far{\aa}s, and I.~Nopens.
\newblock On reliable and unreliable numerical methods for the simulation of
  secondary settling tanks in wastewater treatment.
\newblock {\em Computers Chem. Eng.}, 41:93--105, 2012.

\bibitem{Burger&E&K&L2000}
R.~B{\"u}rger, S.~Evje, K.~H. Karlsen, and K.-A. Lie.
\newblock Numerical methods for the simulation of the settling of flocculated
  suspensions.
\newblock {\em Chem. Eng. J.}, 80:91--104, 2000.

\bibitem{Osores2019}
R.~B\"{u}rger, E.D. Fern{\'{a}}ndez-Nieto, and V.~Osores.
\newblock A dynamic multilayer shallow water model for polydisperse
  sedimentation.
\newblock {\em ESAIM: Math.~Model.~Numer.~Anal.}, 53(4):1391--1432, 2019.

\bibitem{Osores2020}
R.~B\"{u}rger, E.D. Fern{\'{a}}ndez-Nieto, and V.~Osores.
\newblock A multilayer shallow water approach for polydisperse~sedimentation
  with sediment compressibility and mixture viscosity.
\newblock {\em J.~Sci.~Comput.}, 85(2), 2020.

\bibitem{Burger&K&T2005a}
R.~B{\"u}rger, K.~H. Karlsen, and J.~D. Towers.
\newblock A model of continuous sedimentation of flocculated suspensions in
  clarifier-thickener units.
\newblock {\em SIAM J.\ Appl.\ Math.}, 65:882--940, 2005.

\bibitem{Burger2012b}
R.~B\"{u}rger, R.~Ruiz-Baier, K.~Schneider, and H.~Torres.
\newblock A multiresolution method for the simulation of sedimentation in
  inclined channels.
\newblock {\em Int. J. Numer. Anal. Model.}, 9(3):479--504, 2012.

\bibitem{Burger&RB&T2012}
R.~B{\"u}rger, R.~Ruiz-Baier, and H.~Torres.
\newblock A stabilized finite volume element formulation for
  sedimentation-consolidation processes.
\newblock {\em {SIAM} J. Sci. Comput.}, 34:B265--B289, 2012.

\bibitem{Careaga2023}
J.~Careaga and G.N. Gatica.
\newblock Coupled mixed finite element and finite volume methods for a solid
  velocity-based model of multidimensional sedimentation.
\newblock {\em ESAIM: Math.~Model.~Numer.~Anal.}, 57:2529--2556, 2023.

\bibitem{Chancelier1994}
J.-Ph. Chancelier, M.~{Cohen de Lara}, and F.~Pacard.
\newblock Analysis of a conservation {PDE} with discontinuous flux: a model of
  settler.
\newblock {\em SIAM J. Appl. Math.}, 54(4):954--995, 1994.

\bibitem{Castro2012}
M.J.~Castro D{\'{\i}}az and E.~Fern{\'{a}}ndez-Nieto.
\newblock A class of computationally fast first order finite volume solvers:
  {PVM} methods.
\newblock {\em {SIAM} Journal on Scientific Computing}, 34(4):A2173--A2196,
  2012.

\bibitem{CastroDaz2012}
M.J.~Castro D{\'i}az, E.D. Fern{\'{a}}ndez-Nieto, T.~Morales de~Luna,
  G.~Narbona-Reina, and C.~Par{\'{e}}s.
\newblock A {HLLC} scheme for nonconservative hyperbolic problems. application
  to turbidity currents with sediment transport.
\newblock {\em ESAIM: Math.~Model.~Numer.~Anal.}, 47(1):1--32, 2012.

\bibitem{CastroDaz2009}
M.J.~Castro D{\'i}az, E.D. Fern{\'{a}}ndez-Nieto, A.M. Ferreiro, and
  C.~Par{\'{e}}s.
\newblock Two-dimensional sediment transport models in shallow water equations.
  a second order finite volume approach on unstructured meshes.
\newblock {\em Comput.~Methods Appl.~Mech.~Eng.}, 198(33-36):2520--2538, 2009.

\bibitem{SDsiam3}
S.~Diehl.
\newblock Dynamic and steady-state behavior of continuous sedimentation.
\newblock {\em SIAM J. Appl. Math.}, 57(4):991--1018, 1997.

\bibitem{SDwatres2}
S.~Diehl.
\newblock The solids-flux theory -- confirmation and extension by using partial
  differential equations.
\newblock {\em Water Res.}, 42(20):4976--4988, 2008.

\bibitem{Dupont1992}
R.~Dupont and M.~Henze.
\newblock Modelling of the secondary clarifier combined with the activated
  sludge model no.~{\rm 1}.
\newblock {\em Water Sci. Tech.}, 25(6):285--300, 1992.

\bibitem{Ekama1997}
G.A. Ekama, J.L. Barnard, F.W. G{\"u}nthert, P.~Krebs, J.A. McCorquodale, D.S.
  Parker, and E.J. Wahlberg.
\newblock {\em Secondary Settling Tanks: Theory, Modelling, Design and
  Operation}.
\newblock IAWQ scientific and technical report no.~6. International Association
  on Water Quality, England, 1997.

\bibitem{Flores2008}
X.~Flores-Alsina, I.~Rodriguez-Roda, G.~Sin, and K.V. Gernaey.
\newblock Multi-criteria evaluation of wastewater treatment plant control
  strategies under uncertainty.
\newblock {\em Water Res.}, 42(17):4485--4497, 2008.

\bibitem{Gernaey2004}
K.V. Gernaey, M.C.M. van Loosdrecht, M.~Henze, M.~Lind, and S.B.J{\o}rgensen.
\newblock Activated sludge wastewater treatment plant modelling and simulation:
  state of the art.
\newblock {\em Environ. Model. Software}, 19(9):763--783, 2004.

\bibitem{Gustavsson2000}
K.~Gustavsson and J.~Oppelstrup.
\newblock Consolidation of concentrated suspensions {\textendash} numerical
  simulations using a two-phase fluid model.
\newblock {\em Comput. Visual. Sci.}, 3(1--2):39--45, 2000.

\bibitem{Henze1987WR}
M.~Henze, C.P.L. Grady, W.~Gujer, G.V.R. Marais, and T.~Matsuo.
\newblock A general model for single-sludge wastewater treatment systems.
\newblock {\em Water Res.}, 21(5):505--515, 1987.

\bibitem{Henze2000ASMbook}
M.~Henze, W.~Gujer, T.~Mino, and M.C.M. van Loosdrecht.
\newblock {\em Activated Sludge Models {ASM1}, {ASM2}, {ASM2d} and {ASM3}}.
\newblock IWA Publishing, London, UK, 2000.
\newblock IWA Scientific and Technical Report No. 9.

\bibitem{Hu2003}
Z.R. Hu, M.C. Wentzel, and G.A. Ekama.
\newblock {Modelling biological nutrient removal activated sludge systems -- a
  review}.
\newblock {\em Water Res.}, {37}({14}):{3430--3444}, {2003}.

\bibitem{Jayaweera1996}
M.~Jayaweera and T.~Asaeda.
\newblock Modeling of biomanipulation in shallow, eutrophic lakes: An
  application to lake bleiswijkse zoom, the netherlands.
\newblock {\em Ecol.~Model.}, 85(2–3):113–127, 1996.

\bibitem{Jorgensen1988}
S.E. J{\o}rgensen.
\newblock {\em Modelling Eutrophication of Shallow Lakes}, page 177–188.
\newblock Elsevier, 1988.

\bibitem{Li05}
Xiaoye~S. Li.
\newblock An overview of {SuperLU}: Algorithms, implementation, and user
  interface.
\newblock {\em ACM Transactions on Mathematical Software}, 31(3):302--325,
  2005.

\bibitem{Lockett1979}
M.J. Lockett and K.S. Bassoon.
\newblock Sedimentation of binary particle mixtures.
\newblock {\em Powder Technology}, 24(1):1--7, 1979.

\bibitem{Masliyah1979}
J.H. Masliyah.
\newblock Hindered settling in a multi-species particle system.
\newblock {\em Chemical Engineering Science}, 34(9):1166--1168, 1979.

\bibitem{DalMaso1995}
G.~Dal Maso, P.~Le Floch, and F.~Murat.
\newblock Definition and weak stability of nonconservative products.
\newblock {\em J.~Math.~Pures Appl.}, 74:483--548, 1995.

\bibitem{Pauer2000}
J.~Pauer.
\newblock Nitrification in the water column and sediment of a hypereutrophic
  lake and adjoining river system.
\newblock {\em Water Res.}, 34(4):1247–1254, 2000.

\bibitem{Rao2002}
R.~Rao, L.~Mondy, A.~Sun, and S.~Altobelli.
\newblock A numerical and experimental study of batch sedimentation and viscous
  resuspension.
\newblock {\em Int. J. Numer. Methods Fluids.}, 39(6):465--483, 2002.

\bibitem{Tory2003}
E.M. Tory, K.H. Karlsen, R.~B\"{u}rger, and S.~Berres.
\newblock Strongly degenerate parabolic-hyperbolic systems modeling
  polydisperse sedimentation with compression.
\newblock {\em SIAM J.~Appl.~Math.}, 64(1):41--80, 2003.

\bibitem{vanderMolen1994}
D.T. van~der Molen, F.J. Los, L.~van Ballegooijen, and M.P. van~der Vat.
\newblock Mathematical modelling as a tool for management in eutrophication
  control of shallow lakes.
\newblock {\em Hydrobiologia}, 275–276(1):479–492, 1994.

\bibitem{Watts1996a}
R.W. Watts, S.A. Svoronos, and B.~Koopman.
\newblock One-dimensional modeling of secondary clarifiers using a
  concentration and feed velocity-dependent dispersion coefficient.
\newblock {\em Water Res.}, 30(9):2112--2124, 1996.

\bibitem{Wintermeyer2017}
N.~Wintermeyer, A.R. Winters, G.J. Gassner, and D.A. Kopriva.
\newblock An entropy stable nodal discontinuous galerkin method for the two
  dimensional shallow water equations on unstructured curvilinear meshes with
  discontinuous bathymetry.
\newblock {\em J.Comput.~Phys.}, 340:200–242, 2017.

\end{thebibliography}
\end{document}